\documentclass[11pt]{amsart}
\usepackage{graphicx}
\usepackage{amssymb}
\usepackage{amsfonts}
\usepackage[]{amsmath}
\usepackage[]{epsfig}
\usepackage[]{pstricks}
\newpsobject{malla}{psgrid}{subgriddiv=1,griddots=10,gridlabels=6pt}
\usepackage[]{float}
\newpsobject{malla}{psgrid}{subgriddiv=1,griddots=10,gridlabels=6pt}
\newtheorem{theorem}{Theorem}[section]
\newtheorem{lemma}{Lemma}[section]

\newtheorem{remark}{Remark}[section]
\newtheorem{proposition}{Proposition}[section]

\numberwithin{equation}{section}
\newcommand{\R}{\mathbb R}

\newcommand{\p}{\partial}

\begin{document}

\title[ZR system]
{Well-posedness for the 1D Zakharov-Rubenchik system}

\author{Felipe Linares}
\address{Felipe Linares: IMPA, Estrada Dona Castorina 110, Rio de Janeiro, 22460-320, Brazil}
\email{linares@impa.br.}

\author{Carlos Matheus}
\address{ Carlos Matheus: College de France 3, Rue d'Ulm, Paris, France}
\email{cmateus@impa.br}

\begin{abstract} Local and global  well-posedness results are established for the initial value
problem associated to the 1D Zakharov-Rubenchik  system. We show that our results are sharp
in some situations by proving Ill-posedness results otherwise. The global results allow us to study
the norm growth of solutions corresponding to the Schr\"odinger equation term.  We use ideas
recently introduced to study  the classical Zakharov systems.
\end{abstract}

\date{Sept. 5, 2008.}

\maketitle
%%%%%%%%%%%%%%%%%%%%%%%%%%%%%%%%%%%%%%%%%%%%%%%%%%
\section{Introduction}
%%%%%%%%%%%%%%%%%%%%%%%%%%%%%%%%%%%%%%%%%%%%%%%%%%

In this paper we will deal with issues concerning well-posedness for
the initial value problem (IVP) associated to the Zakharov-Rubenchik
system, that is,
\begin{equation}\label{zr1d}
\begin{cases}
i\p_t B+\omega\p^2_x B=\gamma(u-\frac12\nu\,\rho+q\,|B|^2)B,\;\;\;x,\,t\in\R,\\
\theta\,\p_t \rho +\p_x(u-\nu\,\rho)=-\gamma\p_x(|B|^2),\\
\theta\,\p_tu +\p_x(\beta\,\rho-\nu\,u)=\frac12 \gamma\p_x(|B|^2),
\end{cases}
\end{equation}
where $B$ is a complex function, $\rho$ and $u$ are real functions,
$\theta\neq 0$,$\gamma$, $\omega$,  real numbers, $\beta>0$, $\beta-\nu^2\neq 0$, and
$q=\gamma+\nu(\gamma\nu-1)/2(\beta-\nu^2)$. This system is the 1D
version of the most general system deduced by Zakharov and Rubenchik \cite{ZaRu}
to describe the interaction of spectrally narrow high-frequency wave packets of small amplitude
with low-frequency acoustic type oscillations. It has the following form
\begin{equation}\label{zr3d}
\begin{cases}
i\partial_t\psi+iv_g\partial_z\psi=-\frac{\omega''}{2}\,
\partial^2_z\psi-\frac{v_g}{2k}\Delta_\perp\psi+(q|\psi|^2+\beta\rho
+ \alpha\partial_z\varphi)\psi,\\
\partial_t\rho+\rho_0\Delta\varphi+\alpha\partial_z|\psi|^2=0,\\
\partial_t\varphi+\frac{c^2}{\rho_0}\,\rho+\beta|\psi|^2=0,
\end{cases}
\end{equation}
where $\psi$ is a complex function denoting the complex amplitude of the high-frequency carrywave,
$\rho$, $\varphi$ are real functions denoting the density fluctuation and the hydrodynamic potential
respectively. $\alpha$, $\beta$, $q\in \R$, and $\Delta_{\perp}=\p_x^2+\partial_y^2$.

Concerning well-posedness for the IVP associated to \eqref{zr1d},
Oliveira \cite{O} proved local and global well-posedness for data in
$H^2(\R)\times H^1(\R)\times H^1(\R)$. He also studied the existence
and orbital stability of solitary wave solutions for \eqref{zr1d}.
The method used in \cite{O} to establish local well-posedness
follows the ideas of Tsutsumi and Ozawa \cite{ot} to treat the
classical Zakharov systems,
\begin{equation}\label{zsys}
\begin{cases}
i\p_t u+\Delta u=u\,N, \quad\quad x\in\R^n,\,t\in\R,\;\; n\ge 1,\\
\p^2_t N -\Delta N=\Delta(|u|^2).
\end{cases}
\end{equation}

For system \eqref{zr3d}, Ponce and Saut  \cite{ps} proved that the Cauchy problem is locally
well-posed in $H^s(\R^n)$, $\,s>n/2$,  in space dimension $n=2,\,3$. There are several
open questions regarding this system such as global existence of solutions.

In the last few year progress has been made to understand the behavior of solutions
for the Zakharov system \eqref{zsys}.  Here we will apply and extend some of the new
techniques introduced to study system \eqref{zsys} to obtain the results we will describe next.

The first issue we investigate is related to the local well-posedness theory. Before stating
our results in this direction, we make the following change of variables:

Setting
\begin{equation}\label{changevar}
\rho=\psi_1+\psi_2, \text{\hskip10pt and \hskip10pt} u=\sqrt{\beta}\,(\psi_1-\psi_2)
\end{equation}
we write system \eqref{zr1d} as
\begin{equation}\label{zrnew}
\begin{cases}
i\partial_t B + \omega\partial_x^2 B = \gamma\left(\sqrt{\beta}-\frac{\nu}{2}\right)\psi_1 B -
\gamma\left(\sqrt{\beta}+\frac{\nu}{2}\right)\psi_2 B + \gamma q|B|^2 B, \\
\theta\partial_{t}\psi_1 + (\sqrt{\beta}-\nu)\partial_{x}\psi_1 =
\frac{\gamma}{2}\left(-1+\frac{\nu}{2\sqrt{\beta}}\right)\partial_{x}(|B|^2),\\
\theta\partial_t\psi_2 - (\nu+\sqrt{\beta})\partial_x\psi_2 =
\frac{\gamma}{2}\left(-1-\frac{\nu}{2\sqrt{\beta}}\right)\partial_x(|B|^2).
\end{cases}
\end{equation}

Observe that if $(B,\psi_1,\psi_2)$ is a solution of \eqref{zrnew}
with initial data $(B_0,\psi_{10},\psi_{20})$, then
$(B^{\lambda},\psi_1^{\lambda},\psi_2^{\lambda})
=(\lambda\,B,\lambda^2\,\psi_1,\lambda^2\,\psi_2)$ is also a
solution of \eqref{zrnew} with data
$(\lambda\,B_0,\lambda^2\,\psi_{10},\lambda^2\,\psi_{20})$. Hence a
scaling argument suggests local well-posedness for the IVP
\eqref{zrnew} for data in $H^k(\R)\times H^s(\R)\times H^l(\R)$ for
$k>-1/2$, $s,\,l\ge -3/2$.

The local well-posedness theory for \eqref{zrnew} is as follows:

\begin{theorem}\label{t.local}  The Zakharov-Rubenchik
system \eqref{zrnew} is locally well-posed for initial data
\begin{equation*}
(B_0,\psi_{10},\psi_{20})\in H^k(\R)\times H^l(\R)\times H^{s}(\R),
\end{equation*}
where
\begin{equation}\label{e.indices}
\begin{split}
&-\tfrac{1}{2}<k-l\leq 1 , \quad 0\leq l+\tfrac{1}{2}\leq 2k, \\
&-\tfrac{1}{2}<k-s\leq 1 , \quad 0\leq s+\tfrac{1}{2}\leq 2k.
\end{split}
\end{equation}
\end{theorem}

\begin{remark}  Notice that the \lq\lq critical'' indices of the Sobolev spaces where local
well-posed\-ness is expected, i.e., $(s,k,l)=(-1/2,-3/2,-3/2)$ are not reached in our theorem.
\end{remark}

\begin{remark}  Since solutions of system \eqref{zr1d} satisfy the conserved quantities
\begin{align}
I_1(t)&=\int\limits_{\R} |B|^2\,dx,\label{ene1}\\
I_2(t)&=\frac{\omega}2\int\limits_{\R}|B_x|^2+\frac{\gamma q}{4}\int\limits_{\R}|B|^2+\frac{\gamma}2
\int\limits_{\R}(u-\frac{\nu}{2})\rho|B|^2+\frac{\beta}{4}\int\limits_{\R}|\rho|^2+\frac14
\int\limits_{\R}|u|^2-\frac{\nu}2\int\limits_{\R}u\rho,\label{ene2}\\
I_3(t)&=\int\limits_{\R}u\rho\,dx+\frac{i}{2}\int\limits_{\R}(B\bar{B}_x-B_x\bar{B})\,dx,\label{ene3}
\end{align}
and
\begin{equation}\label{ene4}
\begin{split}
I_4(t)=I_2(t)+\tfrac{\nu}{2\theta}I_3(t)= &\,\tfrac{\omega}{2}\int
|B_x|^2 dx +\tfrac{\gamma q}{4}\int |B|^4 dx + \tfrac{\gamma}{2}\int
\left(u-\tfrac{\nu}{2}\rho\right)|B|^2 dx \\&+ \tfrac{\beta}{4}\int
|\rho|^2 dx + \tfrac{1}{4}\int |u|^2 dx + \tfrac{i\nu}{4\theta}\int
(B\overline{B}_x - B_x\overline{B})dx.
\end{split}
\end{equation}
(see \cite{O}) assuming $\omega>0$ and $\beta-\nu^2>0$ we can deduced the global existence of solutions for
data in the space $H^1(\R)\times L^2(\R)\times L^2(\R)$ (see
Proposition \ref{p.globalpoly} below). In particular, this result
implies the stability of solitary wave solutions proved in \cite{O}
for data in the energy space.
\end{remark}

To prove the local result we will follow the scheme used by Ginibre,
Tsutsumi and Velo\cite{GTV} to establish well-posedness for the IVP
associated to the Zakharov system \eqref{zsys}. They used Bourgain
and Kenig, Ponce and Vega arguments.  Since system  \eqref{zrnew} also
contains a cubic term of the function $B$ we need to have
good trilinear estimates in addition to bilinear estimates already
used in \cite{GTV}.

Since solutions of system \eqref{zrnew} satisfy that the $L^2$-norm of $B$ is invariant (\eqref{ene2}), a
natural question regarding global well-posedness arises. Can we extend $B$
to any time?  The answer is positive. Moreover,

\begin{theorem}\label{t.global} The Zakharov-Rubenchik
system \eqref{zrnew} is globally well-posed for initial data
\begin{equation*}
(B_0,\psi_{10},\psi_{20})\in H^k(\R)\times H^l(\R)\times H^{l}(\R),
\end{equation*}
where $0\leq k=l+\tfrac{1}{2}$.
\end{theorem}

\begin{remark} Notice that we establish global existence of solutions for data
in $L^2(\R)\times H^{-1/2}(\R)\times H^{-1/2}(\R)$. Moreover, we can prove
global well-posedness for data with Sobolev indices satisfying $k=l+\frac12$, $k\ge 0$.
In \cite{CHT} global well-posedness for \eqref{zsys} is only presented in the
extremal point $(0,-1/2,-3/2)$.
\end{remark}

To prove Theorem \ref{t.global} we will use the arguments of Colliander, Holmer and Tzirakis in \cite{CHT}
recently put forward to construct global solutions for the 1D Zakharov system \eqref{zsys} in
a similar situation. This is based in the conservation property \eqref{ene1} and the local theory.

Notice that we can write the system \eqref{zrnew} in its integral equivalent form
\begin{equation}\label{isf}
\begin{cases}
B(t)=U(t)(B_0) + i \int\limits_0^t U(t-t')\,\big(|B|^2 + \psi_1+\psi_2\big)(t')\, B(t')\, dt',\\
\psi_{1}(t)=W_+(t)(\psi_{10})+\int\limits_0^t W_+(t-t')\partial_x(|B|^2)(t')\,dt',\\
\psi_{2}(t)=W_-(t)(\psi_{10})-\int\limits_0^t W_-(t-t')\partial_x(|B|^2)(t')\,dt'.
\end{cases}
\end{equation}

The idea of the proof is to perform an iteration scheme. We describe the iteration process
next only considering the extremal case $(0,-1/2,-1/2)$. One of the key observation is that
the interaction of the functions $\psi_1$ and $\psi_2$ in the second and third equations are only with
the function $B$. In other words, we can rewrite \eqref{isf} as
\begin{equation}\label{single}
B(t)=U(t)(B_0) + i \int\limits_0^t U(t-t')\,\big(|B|^2 + \psi_1+\psi_2\big)(t')\, B(t') dt',
\end{equation}
where $\psi_{1}(t)$ and $\psi_{1}(t)$ are as in \eqref{isf}. So one will try to
control the growth of the $L^2$-norm of $B$ using the conserved quantity \eqref{ene1}
and controlling the
growths of  $\psi_1$ and $\psi_2$ in the corresponding $H^l$-norms.
The iteration scheme is as follows given the time of local
well-posedness $T_1$ we denote $B(T_1)=B_1$, $\psi_1(T_1)=\psi_{11}$
and $\psi_1(T_1)=\psi_{21} $. Find $T_2$ via the local result and
then iterate the local theory until time $T_j+T_{j-1}+\dots+T_1$.
These $T_j$ may shrink due to the growth of $\|B_j\|_{L^2}$ and
$\|\psi_{kj}\|_{H^{-1/2}}$, $k=1,2$. To iterate we have to remake the
local theory and then use some spaces with less regularity to
perform the iteration. Since $\|B(t)\|_{L^2}=\|B_0\|_{L^2}$ for all
$t$. The reduction of $T_j$ is then only forced through the growth
of $\|\psi_{kj}\|_{H^{-1/2}}$, $k=1,2$. We then consider \eqref{isf}
posed at $t=T_j$, $(B, \psi_1,\psi_2)(0)=(B_j,\psi_{1j},\psi_{2j})$. Here we may have two situations:
\begin{enumerate}
\item $\|\psi_{1j}\|_{H^{-1/2}}\gg \|B_j\|_{L^2}$ and $\|\psi_{2j}\|_{H^{-1/2}}\gg\|B_j\|_{L^2}$ or
\item $\min \{\|\psi_{1j}\|_{H^{-1/2}},\;\|\psi_{2j}\|_{H^{-1/2}}\}\lesssim \|B_j\|_{L^{2}}$.
\end{enumerate}
We will restrict to discuss possibility (1). On $[T_j, T_{j+1}]$ we write $\psi_{1}$ and $\psi_{2}$
in terms of $\psi_{10}$, $\psi_{20}$ and $B$
\begin{equation}\label{subsys}
\begin{cases}
\psi_{1}(t)=W_+(t-T_j)(\psi_{10})+\int\limits_0^t W_+(t-t')\partial_x(|B|^2)(t')\,dt',\\
\psi_{2}(t)=W_-(t-T_j)(\psi_{20})-\int\limits_0^t W_-(t-t')\partial_x(|B|^2)(t')\,dt'.
\end{cases}
\end{equation}
As $t$ moves from $t=T_j$ to $t=T_{j+1}$, the first terms in both equations in \eqref{subsys}
stay same size in the $H^l$-norm. Any growth in $H^k$-norm as $t$ evolves thus arises from
the second terms. But these terms are intuitively controlled by the conserved size of $B$
in the $L^2$-norm. Also these terms should be small if $[T_j,T_{j+1}]$ is small since
it takes a while for $B$ to contribute to size growth of $\psi_1$ and $\psi_2$. Indeed, the
following estimates hold,
\begin{equation}\label{epsi}
\begin{split}
\|\psi_1\|_{H^{-1/2}}\le c\|\psi_{1j}\|_{H^{-1/2}}+ c\, \Delta T^{\epsilon} \|B_0\|_{L^2}^2\\
\|\psi_2\|_{H^{-1/2}}\le c\|\psi_{2j}\|_{H^{-1/2}}+ c\, \Delta T^{\epsilon} \|B_0\|_{L^2}^2.
\end{split}
\end{equation}
We then iterate the local well-posedness norm with uniform steps of size
\begin{equation}
\Delta T\sim \min\{\|\psi_{1j}\|_{H^{-1/2}}^{-\alpha}, \|\psi_{2j}\|_{H^{-1/2}}^{-\alpha}\}
\end{equation}
over $m$ steps with
\begin{equation}
m\sim \frac{\min\{\|\psi_{1j}\|_{H^{-1/2}}, \|\psi_{2j}\|_{H^{-1/2}}\}}{T^{\epsilon}\|B_0\|_{L^2}^2}.
\end{equation}
We thus extend the solution past $T_j$ to $[T_j, T_j+ m \Delta T]$ without doubling the size
of $\|\psi_{kj}\|_{H^{-1/2}}$, $k=1,2$. A simple computation
reveals
\begin{equation}
m\,\Delta T \sim \|B_0\|^{-2}\, \min\{\|\psi_{10}\|, \|\psi_{20}\|\}^{1-\alpha+\epsilon\alpha}.
\end{equation}
If $1-\alpha+\epsilon\alpha\ge 0$, we obtain global well-posedness by iterating the whole
process. We notice that the factor $\epsilon$ in \eqref{epsi} can be obtained via a new
local theory where the parameters $b,c$ in the $X^{s,b}$ and $W_{\pm}^{l,c}$ space are
not necessarily greater than $1/2$. But we still can show that \eqref{ene1} holds.

The next issue that we are concerned is the study the growth of the $H^s$-norm for solutions of \eqref{zrnew},
corresponding to data in $H^s(\R)$ for noninteger values of $s$. The presence of the conserved quantities
\eqref{ene1}-\eqref{ene4}
and the local existence theory allow us to obtain upper \lq\lq polynomial'' bounds, for the
$H^s$-norm of these solutions.

Energy type estimates were previously used to show that, for solutions of the IVP associated to the
nonlinear Schr\"odinger equation, the $H^s$-norm of these solutions
have an exponential bound, i.e.,
\begin{equation}\label{growth}
\sup_{t\in[0,T]}\|u(t)\|_{H^s}\le c^{|t|}.
\end{equation}
This can be deduced right away from the local existence theory since
\begin{equation}\label{polgrow}
\sup_{[t\in[0,T]}\|u(t)\|_{H^s}\le \|\phi\|_{H^s}+c\,\|\phi\|_{H^s}.
\end{equation}

Bourgain \cite{bou} proved polynomial bounds for certain Hamiltonian PDE's, including the
nonlinear Schr\"odinger equation and the generalized KdV in the periodic setting. He observed
that a slight improvement of \eqref{polgrow}
\begin{equation}
\sup_{[t\in[0,T]}\|u(t)\|_{H^s}\le \|\phi\|_{H^s}+c\,\|\phi\|_{H^s}^{1-\delta},\;\;0<\delta<1,
\end{equation}
implies the polynomial bound $\|u(t)\|_{H^s}\le c\,|t|^{1/\delta}$.
Staffilani \cite{S} showed polynomial bounds for the nonlinear Schr\"odinger and the
generalized on the line. In \cite{CS} Colliander and Staffilani  established
polynomial bounds for solutions of the 1D Zakharov system \eqref{zsys}. Following
their arguments we prove the following:

\begin{proposition}\label{p.polybound} For initial data $B_0,\psi_{10},\psi_{20}\in\mathcal{S}$, the global
solution of \eqref{zrnew} satisfies
\begin{equation}\label{g-sobol}
\|B(t)\|_{H^s}\lesssim 1+|t|^{(s-1)+}.
\end{equation}
\end{proposition}

Finally, we would like to know whether the local well-posedness results obtained here are sharp.
In this direction, we should present some ill-posedness results for the IVP associated to
\eqref{zrnew}. Our results are inspired for those recently obtained
by Holmer \cite{H} for the 1D Zakharov system \eqref{zsys}.

Our first result guarantees that the local result in Theorem \ref{t.local} are the best
possible when $0<k<1$, $l>2k-\tfrac{1}{2}$ or $k\leq 0$, $l>-1/2$. To do this we use the notion introduced
by Christ, Colliander and Tao \cite{cct} called norm-inflation. More precisely, we have the next theorem.

\begin{theorem}\label{t.ill-inflation}Assume that either $0<k<1$, $l>2k-\tfrac{1}{2}$ or $k\leq 0$, $l>-1/2$.
Then, there are some constants $T>0$, $\alpha=\alpha(k,s)$ and a sequence $f_N\in\mathcal{S}$ with
$\|f_N\|_{H^k}\leq 1$ for all $N\in\mathbb{N}$ such that the solution $(B_N,\phi_{1,N},\phi_{2,N})$ of the
Zakharov-Rubenchik system~\eqref{zrnew} on $[0,T]$ with initial data $(f_N,0,0)$ verifies
$$
\|\psi_{1,N}(t)\|_{H^l}\gtrsim t\cdot N^\alpha
$$
for $0<t\leq T$ and $N\gtrsim \tfrac{1}{t}$.

Similarly, if either $0<k<1$, $s>2k-\tfrac{1}{2}$ or $k\leq 0$, $s>-1/2$, then we can find some constants
$T>0$, $\alpha=\alpha(k,s)$ and a sequence $g_N\in\mathcal{S}$ with $\|g_N\|_{H^k}\leq 1$ for all $N$ such
that the solution $(B_N,\phi_{1,N},\phi_{2,N})$ of~\eqref{zrnew} on $[0,T]$ with initial data $(g_N,0,0)$
verifies
$$
\|\psi_{2,N}(t)\|\gtrsim t\cdot N^{\alpha}
$$
for $t\in (0,T]$ and $N\gtrsim\tfrac{1}{t}$.
\end{theorem}

The second result is weaker than the previous one. It only affirms that a local result cannot be obtained
by Picard iteration. This was first proved by Bourgain for the Korteweg-de Vries equation and Nonlinear
Schr\"odinger \cite{Bo7}. Several improvements have been obtained for other models. In our case, we show
that the best local result suggested for the scaling argument cannot be attained.

\begin{theorem}\label{t.ill-C2} For any $k\in\mathbb{R}$, $l,s$ with $\min\{l,s\}<-1/2$ and $T>0$, the map
data-to-solution associated to the Zakharov-Rubenchik system~\eqref{zrnew} from
$H^k(\R)\times H^l(\R)\times H^s(\R)$
to $C([0,T]; H^k(\R)\times H^l(\R)\times H^s(\R))$ is not $C^2$ at the origin $(0,0,0)$.
\end{theorem}

Finally, we show ill-posedness for data in $H^k(\R)\times H^l(\R)\times H^s(\R))$, $k=0$ and $\min\{l,s\}<-3/2$,
this assures that below the indices suggested by the scaling argument the IVP is in fact ill-posed. Results in this
direction were introduced first by Birnir, Kenig, Ponce, Svanstedt, and Vega \cite{BKPSV}, and generalizations by
Christ, Colliander and Tao \cite{cct}. The result reads as follows.

\begin{theorem}\label{t.ill-decoherence}Assume that $k=0$ and $\min\{l,s\}<-3/2$. Then, given $T>0$ and
$\delta>0$, there exists a pair $(B_0,\psi_{+0},\psi_{-0})$ and $(B_0',\psi_{+0}',\psi_{-0}')$ of initial
data with
\begin{equation*}
\|B_0\|_{H^k},\|\psi_{+0}\|_{H^l},\|\psi_{-0}\|_{H^s},
\|B_0'\|_{H^k},\|\psi_{+0}'\|_{H^l},\|\psi_{-0}'\|_{H^s}\leq 1
\end{equation*}
such that the associated solutions $(B,\psi_+,\psi_-)$ and $(B',\psi_+',\psi_-')$ of the Zakharov-Rubenchik
system~\eqref{zrnew} on the time interval $[0,T]$ start close:
$$\|B_0-B_0'\|_{H^k},\|\psi_{+0}-\psi_{+0}'\|_{H^l}+\|\psi_{-0}-\psi_{-0}'\|_{H^s}\leq \delta$$
but they become separated by time $T$ (in the Schr\"odinger variable):
$$\|B-B'\|_{L_{[0,T]}^{\infty}H_x^k}\sim 1.$$
\end{theorem}

This paper is organized as follows: in Section 2,  we will deal with the local and global theory
for system \eqref{zrnew}. In Section 3, the growth rate of the $H^s$-norm will be established.
Finally, we will show the sharpness of some of the by establishing some ill-posedness
results in Section 4. \\

{\bf Acknowledgements.} We would like to thank German Fonseca for fruitful conversations
regarding this work. We are also grateful to Jim Colliander for some useful comments on
a previous version.

%%%%%%%%%%%%%%%%%%%%%%%%%%%%%%%%%%%%%%%%%%%%%%%%
\section{Well-posedness of Zakharov-Rubenchik system}
%%%%%%%%%%%%%%%%%%%%%%%%%%%%%%%%%%%%%%%%%%%%%%%%
This section is devoted to the proof of the following results Theorem \ref {t.local}
Theorem \ref{t.global}.

We split the proof of these theorems into two steps: firstly, we
verify that the linear and multilinear estimates obtained by
Ginibre, Tsutsumi and Velo~\cite{GTV} directly implies that the
Zakharov-Rubenchik system \eqref{zr1d} is locally well-posed
(Theorem~\ref{t.local}); secondly, we modify this local
well-posedness result and apply the conservation of the $L^2$-mass
of  $B$ (along the lines of a recent work~\cite{CHT} of Colliander,
Holmer and Tzirakis)  to derive the Theorem~\ref{t.global}.

\subsection{Local well-posedness}

We start with the linear and multilinear estimates derived
in~\cite{GTV} for the Zakharov and Benney systems. Let $U(t) =
e^{it\Delta}$ be the free Schr\"odinger linear group and
$W_{\pm}(t)=e^{\pm t\partial_x}$ be the free linear group of the
transport equations.

In the sequel, we use the following norms:
\begin{equation*}
\|f\|_{H^s}:=\left(\int \langle\xi\rangle^{2s} |\widehat{f}(\xi)|^2
d\xi\right)^{1/2},
\end{equation*}
\begin{equation*}
\|u\|_{X^{k,b}}:=\left(\int \langle\xi\rangle^{2k}
\langle\tau+\xi^2\rangle^{2b} |\widehat{u}(\xi,\tau)|^2 d\xi
d\tau\right)^{1/2},
\end{equation*}
\begin{equation*}
\|v\|_{W_{\pm}^{l,c}}:=\left(\int \langle\xi\rangle^{2l}
\langle\tau\pm\xi\rangle^{2c} |\widehat{v}(\xi,\tau)|^2 d\xi
d\tau\right)^{1/2},
\end{equation*}
where $\langle x\rangle:= 1+ |x|$.

Let $\kappa$ be a smooth bump function so that $\kappa(t)=1$ for
$|t|\leq 1$ and $\kappa(t)=0$ for $|t|\geq 2$. Define $\kappa_T(t):=
\kappa(t/T)$. In this setting, Ginibre, Tsutsutmi and
Velo~\cite{GTV} proved the following linear and multilinear
estimates for the Zakharov and Benney systems:

\begin{lemma}[Linear group estimates]\label{l.linear-1}It holds
\begin{itemize}
\item $\|\kappa(t)U(t)u_0\|_{X^{k,b}}\lesssim_b \|u_0\|_{H^k}$ and
\item
$\|\kappa(t)W_{\pm}(t)v_0\|_{W_{\pm}^{l,c}}\lesssim_c \|v_0\|_{H^l}$
\end{itemize}
Moreover, if $0<T<1$ and $-\tfrac{1}{2}< b'\leq b < \tfrac{1}{2}$,
then
\begin{itemize}
\item $\|\kappa_T(t)u\|_{X^{k,b'}}\lesssim_{b',b} T^{b-b'}
\|u\|_{X^{k,b}}$,
\item $\|\kappa_T(t)v\|_{W_{\pm}^{l,b'}}\lesssim_{b',b} T^{b-b'}
\|v\|_{W_{\pm}^{l,b}}$.
\end{itemize}
\end{lemma}

In the sequel, we define $U *_R F(t,x):=\int_0^t U(t-t')F(t',x) dt'$ and $W_{\pm} *_R
F(t,x):= \int_0^t W_{\pm}(t-t')F(t',x) dt'$.

\begin{lemma}[Duhamel estimates]\label{l.linear-2}Let $-\tfrac{1}{2}< -a \leq 0\leq b\leq
1-a$. It holds
\begin{itemize}
\item $\|\kappa_T(t)U *_R F\|_{X^{k,b}}\lesssim T^{1-b-a}
\|F\|_{X^{k,-a}}$,
\item $\|\kappa_T(t)W_{\pm} *_R G\|_{W_{\pm}^{l,b}}\lesssim
T^{1-b-a} \|G\|_{W_{\pm}^{l,-a}}$.
\end{itemize}
\end{lemma}

Next, we recall the following multilinear estimates derived by
Bourgain~\cite{B} and Ginibre, Tsutsumi and Velo~\cite{GTV}:
\begin{lemma}[Multilinear estimates]\label{l.multilinear}It holds
\begin{itemize}
\item for any $k\geq 0$ and $1/2<b<5/8$,
\begin{equation*}
\begin{split}
\|uv\overline{z}\|_{X^{k,b-1}}\lesssim
&\|u\|_{X^{k,3/8+}}\|v\|_{X^{0,3/8+}}\|z\|_{X^{0,3/8+}} +
\|u\|_{X^{0,3/8+}}\|v\|_{X^{k,3/8+}}\|z\|_{X^{0,3/8+}} + \\
&\|u\|_{X^{0,3/8+}}\|v\|_{X^{0,3/8+}}\|z\|_{X^{k,3/8+}}
\end{split}
\end{equation*}
\item for any $k\geq 0$, $l\geq -\tfrac{1}{2}$, $k-l\leq \min\{1,2c\}$, $b,b',c>\tfrac{3}{8}$,
\begin{equation*}
\|\psi_{\pm}u\|_{X^{k,-c}}\lesssim \|\psi_{\pm}\|_{W_{\pm}^{l,b'}}
\|u\|_{X^{k,b}}
\end{equation*}
\item for any $k\geq 0$, $k-l>-1/2$, $2k-l\geq 1/2$,
$b,c>\tfrac{3}{8}$, $l-k+1\leq 2b$, $l+1-k<2c+1/2$,
\begin{equation*}
\|\partial_x(u\overline{v})\|_{W_{\pm}^{l,-c}}\lesssim
\|u\|_{X^{k,b}} \|v\|_{X^{k,b}}.
\end{equation*}
\end{itemize}
\end{lemma}

Combining these lemmas, it is a standard matter to show the local
well-posedness result of Theorem~\ref{t.local}:

\begin{proof}[Proof of Theorem~\ref{t.local}]
We already know that it suffices to consider the
equation \eqref{zrnew}. The integral formulation of \eqref{zrnew} is
\begin{equation*}
\begin{cases}
B(t) = U(t)B_0 +i U *_R \left\{\psi_1 B + \psi_2 B + |B|^2 B\right\}(t), \\
\psi_1(t) = W_{+}(t)(\psi_{10}) + W_{+} *_R
\partial_x(|B|^2)(t),
\\
\psi_2(t) = W_{-}(t)\psi_{20} - W_{-} *_R \partial_x(|B|^2)(t).
\end{cases}
\end{equation*}
Fix $0<T<1$ and define
\begin{equation}\label{e.zr2}
\begin{split}
\Lambda_T (B) (t) := &\kappa(t) U(t) B_0 + i\kappa_T(t) U *_R
\left\{(\kappa W_{+}(\psi_{10}) + W_+ *_R
\partial_x(|B|^2))B\right\}(t) + \\
&i\kappa_T(t) U *_R \left\{(\kappa W_-(\psi_{20}) -
W_- *_R \partial_x(|B|^2)) B\right\}(t) + \\
&i\kappa_T (t) U *_R (|B|^2 B)(t),
\end{split}
\end{equation}
Our task is reduced to find a fixed point $B = \Lambda_T(B)$ of
$\Lambda_T$. Applying the linear and multilinear estimates of
Lemmas~\ref{l.linear-1},~\ref{l.linear-2} and~\ref{l.multilinear},
we get
\begin{equation*}
\|\Lambda_T(B)\|_{X^{k,1/2+}}\lesssim \|B_0\|_{H^k} + T^{3/8-}
\left\{\|\psi_{10}\|_{H^l} + \|\psi_{20}\|_{H^s} +
\|B\|_{X^{k,1/2+}}^2\right\}\cdot\|B\|_{X^{k,1/2+}}
\end{equation*}
and
\begin{equation*}
\begin{split}
&\|\Lambda_T(B)-\Lambda_T(\widetilde{B})\|_{X^{k,1/2+}}\lesssim \\
&T^{3/8-} \left\{\|\psi_{10}\|_{H^l}+\|\psi_{20}\|_{H^s}+
(1+\|B\|_{X^{k,1/2+}})(1+\|\widetilde{B}\|_{X^{k,1/2+}})
\right\}\|B-\widetilde{B}\|_{X^{k,1/2+}}
\end{split}
\end{equation*}
This implies that $\Lambda_T$ is a contraction of a large ball of
$X^{k,1/2+}$ when $T$ is sufficiently small (depending only on
$\|u_0\|_{H^k}, \|\psi_{10}\|_{H^l}, \|\psi_{20}\|_{H^s}$). In
particular, \eqref{zrnew} is locally well-posed for some time
interval $[0,T]$ with $T=T(\|u_0\|_{H^k}, \|\psi_{10}\|_{H^l},
\|\psi_{20}\|_{H^s})>0$.
\end{proof}

\subsection{Global well-posedness} Following the lines of Colliander, Holmer and
Tzirakis~\cite{CHT}, we begin with some refinements of the linear
and multilinear estimates of
Lemmas~\ref{l.linear-1},~\ref{l.linear-2} and~\ref{l.multilinear}:

\begin{lemma}[Linear estimates]\label{l.linear}
For $0<T\leq 1$, $t\in\mathbb{R}$ and $0\leq b,b_1\leq 1/2$, it
holds:
\begin{itemize}
\item $\|U(t)u_0\|_{H^s}=\|u_0\|_{H^s}$ and
$\|\kappa_T(t)U(t)u_0\|_{X^{s,b}}\lesssim
T^{\tfrac{1}{2}-b}\|u_0\|_{H^s}$;
\item $\|W_{\pm}(t)v_0\|_{H^l}=\|v_0\|_{H^l}$ and
$\|\kappa_T(t)W_{\pm}(t)v_0\|_{W_{\pm}^{l,b_1}}\lesssim
T^{\tfrac{1}{2}-b}\|v_0\|_{H^l}$;
\end{itemize}
\end{lemma}
\begin{proof}See Lemma 2.1 of~\cite{CHT}.
\end{proof}

\begin{lemma}[Duhamel estimates]\label{l.linearII}For $0<T\leq 1$, $0\leq c,c_1 < 1/2$,
$0\leq b\leq b+c\leq 1$ and $0\leq b_1\leq b_1+c_1\leq 1$, we have
\begin{itemize}
\item $\|U*_R F\|_{C^0([0,T],H^k)}\lesssim
T^{\tfrac{1}{2}-c}\|F\|_{X^{k,-c}}$ and $\|\kappa_T U*_R
F\|_{X^{k,b}}\lesssim T^{1-b-c}\|F\|_{X^{k,-c}}$;
\item $\|W_{\pm}*_R G\|_{C^0([0,T],H^l)}\lesssim
T^{\tfrac{1}{2}-c_1}\|G\|_{W_{\pm}^{l,-c_1}}$ and $\|W_{\pm}*_R
G\|_{W_{\pm}^{l,b_1}}\lesssim
T^{1-b_1-c_1}\|G\|_{W_{\pm}^{l,-c_1}}$.
\end{itemize}
\end{lemma}
\begin{proof}See Lemma 2.3 of~\cite{CHT}.
\end{proof}

\begin{lemma}[Multilinear estimates]\label{l.multilinear-2}
It holds
\begin{itemize}
\item a) for any $k\geq 0$ and $3/8<c<1/2$,
\begin{equation*}
\begin{split}
\|uv\overline{z}\|_{X^{k,-c}}\lesssim
&\|u\|_{X^{k,3/8+}}\|v\|_{X^{0,3/8+}}\|z\|_{X^{0,3/8+}} +
\|u\|_{X^{0,3/8+}}\|v\|_{X^{k,3/8+}}\|z\|_{X^{0,3/8+}} + \\
&\|u\|_{X^{0,3/8+}}\|v\|_{X^{0,3/8+}}\|z\|_{X^{k,3/8+}}
\end{split}
\end{equation*}
\item b) for any $k\geq 0$, $k-l\leq 1/2$,
$\tfrac{1}{4}<b,c,b_1<\tfrac{1}{2}$ and $b+b_1+c\geq 1$,
\begin{equation*}
\|\psi_{\pm}u\|_{X^{k,-c}}\lesssim \|\psi_{\pm}\|_{W_{\pm}^{l,b_1}}
\|u\|_{X^{0,b}} + \|\psi_{\pm}\|_{W_{\pm}^{-1/2,b_1}}
\|u\|_{X^{k,b}}
\end{equation*}
\item c) for any $l\geq -1/2$, $k-l\geq 1/2$,
$\tfrac{1}{4}<b,c_1<\tfrac{1}{2}$ and $2b+c_1\geq 1$,
\begin{equation*}
\|\partial_x(u\overline{v})\|_{W_{\pm}^{l,-c_1}}\lesssim
\|u\|_{X^{k,b}} \|v\|_{X^{0,b}} + \|u\|_{X^{0,b}} \|v\|_{X^{k,b}}.
\end{equation*}
\end{itemize}
\end{lemma}

\begin{proof}From Lemma~\ref{l.multilinear}, we already
know that a) holds. Furthermore, from lemma 3.1 of Colliander,
Holmer and Tzirakis~\cite{CHT}, it holds
\begin{itemize}
\item $\|\psi_{\pm}u\|_{X^{0,-c}}\lesssim
\|\psi_{\pm}\|_{W_{\pm}^{-1/2,b_1}}\|u\|_{X^{0,b}}$ if
$\tfrac{1}{4}<b,c,b_1<\tfrac{1}{2}$ and $b+b_1+c\geq 1$;
\item $\|\partial_x(u_1\overline{u_2})\|_{W_{\pm}^{-1/2,-c_1}}\lesssim
\|u_1\|_{X^{0,b}} \|u_2\|_{X^{0,b}}$ if
$\tfrac{1}{4}<b,c_1<\tfrac{1}{2}$ and $2b+c_1\geq 1$.
\end{itemize}
Also, the triangular inequality implies $\langle\xi\rangle^{r}\leq
\langle\xi\rangle^{r'} \langle\xi_1\rangle^{r-r'} +
\langle\xi\rangle^{r'} \langle\xi_2\rangle^{r-r'}$, where $r\geq r'$
and $\xi=\xi_1+\xi_2$. Denoting by $\widehat{J^s
u}(\xi):=\langle\xi\rangle^s \widehat{u}(\xi)$ and putting together
these facts, we obtain
\begin{equation*}
\begin{split}
\|\psi_{\pm}u\|_{X^{k,-c}}&\lesssim \|u\cdot J^k\psi_{\pm}\|_{X^{0,-c}}
+ \|\psi_{\pm}\cdot J^k u\|_{X^{0,-c}} \\
&\lesssim \|J^k\psi_{\pm}\|_{W_{\pm}^{-1/2,b_1}} \|u\|_{X^{0,b}} +
\|\psi_{\pm}\|_{W_{\pm}^{-1/2,b_1}} \|J^k u\|_{X^{0,b}} \\
&\lesssim \|\psi_{\pm}\|_{W_{\pm}^{l,b_1}} \|u\|_{X^{0,b}} +
\|\psi_{\pm}\|_{W_{\pm}^{-1/2,b_1}} \|u\|_{X^{k,b}},
\end{split}
\end{equation*}
if $k\geq 0$, $k-l\leq 1/2$, $\tfrac{1}{4}<b,c,b_1<\tfrac{1}{2}$,
$b+b_1+c\geq 1$, and
\begin{equation*}
\begin{split}
\|\partial_x(u\overline{v})\|_{W_{\pm}^{l,-c}}&\lesssim
\|\partial_x(J^{l+1/2}u\cdot \overline{v})\|_{W_{\pm}^{-1/2,-c}} +
\|\partial_x(u\cdot J^{l+1/2}\overline{v})\|_{W_{\pm}^{-1/2,-c}}\\
&\lesssim \|J^{l+1/2}u\|_{X^{0,b}} \|v\|_{X^{0,b}} +
\|u\|_{X^{0,b}} \|J^{l+1/2}v\|_{X^{0,b}} \\
&\lesssim \|u\|_{X^{k,b}} \|v\|_{X^{0,b}} + \|u\|_{X^{0,b}}
\|v\|_{X^{k,b}},
\end{split}
\end{equation*}
if $k\geq 0$, $k-l\geq 1/2$, $\tfrac{1}{4}<b,c_1<\tfrac{1}{2}$,
$2b+c_1\geq 1$. This completes the proof.
\end{proof}

Next, we recall that the $L^2$ mass of $B$ is a conserved quantity, i.e.,
\begin{equation*}
I_1(t):=\int|B(t)|^2 dx=\int |B_0|^2 dx
\end{equation*}

Now, we combine the previous lemmas and remarks to conclude the proof of
Theorem~\ref{t.global}:

\begin{proof}[Proof of Theorem~\ref{t.global}] Firstly, we treat system \eqref{zrnew} for
initial data $(B_0,\psi_{10},\psi_{20})$ in the space $L^2(\R)\times H^{-1/2}(\R)\times
H^{-1/2}(\R)$. Fix $0<T<1$ and define
\begin{equation}\label{e.global.1}
\begin{split}
&\Phi_T(B,\psi_1,\psi_2)(t):=\kappa_T(t)U(t)B_0 + \kappa_T(t)
U*_R\{(\psi_1+\psi_2+|B|^2)B\}(t) \\
&\Psi_T^{+}(B)(t):=\kappa_T(t) W_+(t)(\psi_{10})+ \kappa_T(t)
W_+*_R\{\partial_x(|B|^2)\}(t) \\
&\Psi_T^-(B)(t):=\kappa_T(t)W_-(t)(\psi_{20}) -
\kappa_T(t)W_-*_R\{\partial_x(|B|^2)\}(t)
\end{split}
\end{equation}
We seek a fixed point
$(B,\psi_1,\psi_2)=(\Phi_T(B,\psi_1,\psi_2),\Psi_T^+(B),\Psi_T^-(B))$.
To do so, we estimate~(\ref{e.global.1}) in $X^{0,b}\times
W_+^{-1/2,b_1}\times W_-^{-1/2,b_1}$ (with $b=c=3/8+\varepsilon$,
$b_1=c_1=1/4+\varepsilon$) using the
Lemmas~\ref{l.linear},~\ref{l.linearII},~\ref{l.multilinear-2}, so
that
\begin{equation*}
\begin{split}
&\|\Phi_T(B,\psi_1,\psi_2)\|_{X^{0,b}}\lesssim
T^{\tfrac{1}{8}-\varepsilon}\|B_0\|_{L^2} \!+\!
T^{\tfrac{1}{4}-2\varepsilon}(\|\psi_1\|_{W_+^{-1/2,b_1}}\!+\!
\|\psi_2\|_{W_-^{-1/2,b_1}}+\|B\|_{X^{0,b}}^2)\|B\|_{X^{0,b}} \\
&\|\Psi_T^+(B)\|_{W_+^{-1/2,b_1}}\lesssim
T^{\tfrac{1}{4}-\varepsilon}\|\psi_{10}\|_{H^{-1/2}} +
T^{\tfrac{1}{2}-2\varepsilon}\|B\|_{X^{0,b}}^2 \\
&\|\Psi_T^-(B)\|_{W_-^{-1/2,b_1}}\lesssim
T^{\tfrac{1}{4}-\varepsilon}\|\psi_{20}\|_{H^{-1/2}} +
T^{\tfrac{1}{2}-2\varepsilon}\|B\|_{X^{0,b}}^2
\end{split}
\end{equation*}
and
\begin{equation*}
\begin{split}
&\|\Phi_T(B,\psi_1,\psi_2)-\Phi_T(\widetilde{B},\widetilde{\psi_1},\widetilde{\psi_2})\|_{X^{0,b}}\lesssim
T^{\tfrac{1}{4}-2\varepsilon}(\|\psi_1\|_{W_+^{-1/2,b_1}}+\|\psi_2\|_{W_-^{-1/2,b_1}})
\|B-\widetilde{B}\|_{X^{0,b}} \\
&+T^{\tfrac{1}{4}-2\varepsilon}\|\widetilde{B}\|_{X^{0,b}}
\|\psi_1-\widetilde{\psi_1}\|_{W_+^{-1/2,b_1}}+T^{\tfrac{1}{4}-2\varepsilon}
\|\widetilde{B}\|_{X^{0,b}}\|\psi_2-\widetilde{\psi_2}\|_{W_-^{-1/2,b_1}}\\
&+T^{\tfrac{1}{4}-2\varepsilon}(\|B\|_{X^{0,b}}+\|\widetilde{B}\|_{X^{0,b}})^2
\|B-\widetilde{B}\|_{X^{0,b}},\\
&\|\Psi_T^+(B)-\Psi_T^+(\widetilde{B})\|_{W_+^{-1/2,b_1}}\lesssim
T^{\tfrac{1}{2}-2\varepsilon}(\|B\|_{X^{0,b}}+\|\widetilde{B}\|_{X^{0,b}})
\|B-\widetilde{B}\|_{X^{0,b}}, \\
&\|\Psi_T^-(B)-\Psi_T^-(\widetilde{B})\|_{W_-^{-1/2,b_1}}\lesssim
T^{\tfrac{1}{2}-2\varepsilon}(\|B\|_{X^{0,b}}+\|\widetilde{B}\|_{X^{0,b}})
\|B-\widetilde{B}\|_{X^{0,b}}.
\end{split}
\end{equation*}
Thus, for any $0<T<1$ such that
\begin{equation*}
\begin{split}
T^{\tfrac{3}{8}-3\varepsilon}\|B_0\|_{L^2}\leq 1, \quad \quad \quad T^{\tfrac{1}{2}-4\varepsilon}
\|B_0\|_{L^2}^2\leq 1,
\end{split}
\end{equation*}
and
\begin{equation}\label{e.global.2}
T^{\tfrac{1}{2}-3\varepsilon}\|\psi_{10}\|_{H^{-1/2}}\leq 1, \quad \quad \quad
T^{\tfrac{1}{2}-3\varepsilon}\|\psi_{20}\|_{H^{-1/2}}\leq 1,
\end{equation}
\begin{equation}\label{e.global.3}
T^{\tfrac{1}{2}-2\varepsilon}\cdot T^{\tfrac{1}{4}-2\varepsilon} \|B_0\|_{L^2}^2 \leq
T^{\tfrac{1}{4}-\varepsilon}\|\psi_{10}\|_{H^{-1/2}}, \quad T^{\tfrac{1}{2}-2\varepsilon}\cdot
T^{\tfrac{1}{4}-2\varepsilon} \|B_0\|_{L^2}^2 \leq T^{\tfrac{1}{4}-\varepsilon}\|\psi_{20}\|_{H^{-1/2}}
\end{equation}
we obtain sufficient conditions for the application of a standard contraction argument yielding a
unique fixed point $(B,\psi_1,\psi_2)\in X^{0,b}\times W_+^{-1/2,b_1}\times W_-^{-1/2,b_1}$ for
(\ref{e.global.1}) verifying
\begin{equation}\label{e.global.4}
\begin{split}
&\|B\|_{X^{0,b}}\lesssim T^{\tfrac{1}{8}-\varepsilon}\|B_0\|_{L^2},\\ &\|\psi_1\|_{W_+^{-1/2,b_1}}\lesssim
T^{\tfrac{1}{4}-\varepsilon}\|\psi_{10}\|_{H^{-1/2}},\\
&\|\psi_2\|_{W_-^{-1/2,b_1}}\lesssim T^{\tfrac{1}{4}-\varepsilon}\|\psi_{20}\|_{H^{-1/2}}.
\end{split}
\end{equation}
Nevertheless, we can estimate $\Phi_T(B,\psi_1,\psi_2)$ in $C^0([0,T],L^2)$ using the
Lemmas~\ref{l.linear},~\ref{l.linearII},~\ref{l.multilinear-2} and~(\ref{e.global.4}) to prove that
$B\in C^0([0,T],L^2)$. By the conservation of the $L^2$ norm of $B$ (see \eqref{ene1}) says
that the norm $\|B(t)\|_{L^2}=\|B_0\|_{L^2}$ is unchanged during the evolution. Therefore, it remains only
to deal with the possible growth of $\|\psi_{10}\|_{H^{-1/2}}$ and $\|\psi_{20}\|_{H^{-1/2}}$. Assume that,
after some iteration, we attain a time $t$ such that either
$\|\psi_1(t)\|_{H^{-1/2}}\gg\|B(t)\|_{L^2}^2=\|B_0\|_{L^2}^2$ or
$\|\psi_2(t)\|_{H^{-1/2}}\gg\|B(t)\|_{L^2}^2=\|B_0\|_{L^2}^2$. In the sequel, the time position $t$ will be
the initial time $t=0$, so that either $\|\psi_{10}\|_{H^{-1/2}}\gg\|B_0\|_{L^2}^2$ or
$\|\psi_{20}\|_{H^{-1/2}}\gg\|B_0\|_{L^2}^2$. We have two possibilities:
\begin{itemize}
\item a) $\|\psi_{10}\|_{H^{-1/2}}\gg\|B_0\|_{L^2}^2$ \textbf{and} $\|\psi_{20}\|_{H^{-1/2}}
\gg\|B_0\|_{L^2}^2$:
it follows that (\ref{e.global.3}) is \emph{automatically} satisfied and we may choose
\begin{equation}\label{e.global.5}
T\sim \min\{\|\psi_{10}\|^{-1/(\tfrac{1}{2}-3\varepsilon)}, \|\psi_{20}\|^{-1/(\tfrac{1}{2}-3\varepsilon)}\}
\end{equation}
in the (local) iteration scheme (by (\ref{e.global.2})). Because
$$\psi_1(t)=W_+(\psi_{10})(t)+W_+*_R(\partial_x(|B|^2)) \quad \textrm{ and }\quad
\psi_2(t)=W_-(\psi_{20})(t)-W_-*_R(\partial_x(|B|^2)),$$
a simple application of the Lemmas \ref{l.linear}, \ref{l.linearII} and (\ref{e.global.4}) yields
$$
\|\psi_1(T)\|_{H^{-1/2}}\leq \|\psi_{10}\|_{H^{-1/2}} + C\cdot T^{\tfrac{1}{2}-3\varepsilon}\|B_0\|_{L^2}^2
$$
and
$$
\|\psi_2(T)\|_{H^{-1/2}}\leq \|\psi_{20}\|_{H^{-1/2}} + C\cdot T^{\tfrac{1}{2}-3\varepsilon}\|B_0\|_{L^2}^2.
$$
These estimates show that we can carry out $m$ iterations on time intervals of size (\ref{e.global.5})
before the quantities $\|\psi_1(t)\|_{H^{-1/2}}$ and $\|\psi_2(t)\|_{H^{-1/2}}$ doubles, where
\begin{equation}\label{e.global.6}
m\sim \frac{\min\{\|\psi_{10}\|_{H^{-1/2}},\|\psi_{20}\|_{H^{-1/2}}\}}
{T^{\tfrac{1}{2}-3\varepsilon}\|B_0\|_{L^2}^2}.
\end{equation}
Observe that, by (\ref{e.global.5}) and (\ref{e.global.6}), the amount of time that we advanced
(after these $m$ iterations) is
$$mT\sim \|B_0\|_{L^2}^{-2}.$$
This proves the global well-posedness result in the case a) and it shows that
$$\|\psi_1(t)\|_{H^{-1/2}}\leq \exp(Ct\|B_0\|_{L^2}^2)\max\{\|\psi_{10}\|_{H^{-1/2}}, \|B_0\|_{L^2}^2\}$$
and
$$\|\psi_2(t)\|_{H^{-1/2}}\leq \exp(Ct\|B_0\|_{L^2}^2)\max\{\|\psi_{20}\|_{H^{-1/2}}, \|B_0\|_{L^2}^2\}.$$
\item b) $\min\{\|\psi_{10}\|_{H^{-1/2}}, \|\psi_{20}\|_{H^{-1/2}}\}\lesssim \|B_0\|_{L^2}^2$: this case is
similar to the previous one, but we need to rework the entire argument above; by interchanging the roles of
$\psi_{10}$ and $\psi_{20}$ if necessary, we can suppose that
$\|\psi_{20}\|_{H^{-1/2}}\lesssim\|B_0\|_{L^2}^2$ and $\|\psi_{10}\|_{H^{-1/2}}\gg\|B_0\|_{L^2}^2$; in this
case, for any $0<T<1$ such that $T^{\tfrac{1}{2}-4\varepsilon}\|B_0\|_{L^2}^2\leq 1$,
\begin{equation}\label{e.global.7}
T^{\tfrac{1}{2}-3\varepsilon}\|\psi_{10}\|_{H^{-1/2}}\leq 1,
\end{equation}
and
\begin{equation}\label{e.global.8}
T^{\tfrac{1}{2}-2\varepsilon}\cdot T^{\tfrac{1}{4}-2\varepsilon} \|B_0\|_{L^2}^2 \leq
T^{\tfrac{1}{4}-\varepsilon}\|\psi_{10}\|_{H^{-1/2}}
\end{equation}
we get a local-in-time solution $(B,\psi_1,\psi_2)$ of (\ref{e.global.1}) on the time interval $T$ satisfying
\begin{equation}\label{e.global.9}
\begin{split}
&\|B\|_{X^{0,b}}\lesssim T^{\tfrac{1}{8}-\varepsilon}\|B_0\|_{L^2},\\ &\|\psi_1\|_{W_+^{-1/2,b_1}}\lesssim
T^{\tfrac{1}{4}-\varepsilon}\|\psi_{10}\|_{H^{-1/2}},\\
&\|\psi_2\|_{W_-^{-1/2,b_1}}\lesssim T^{\tfrac{1}{4}-\varepsilon}\|B_0\|_{L^2}^2.
\end{split}
\end{equation}
Note that (\ref{e.global.8}) is fulfilled because $\|\psi_{10}\|_{H^{-1/2}}\gg\|B_0\|_{L^2}^2$ (by hypothesis).
Hence, there exists a solution of the Zakharov-Rubenchik system on a time interval $[0,T]$ of size
\begin{equation}\label{e.global.10}
T\sim \|\psi_{10}\|^{-1/(\tfrac{1}{2}-3\varepsilon)}
\end{equation}
Again, using the Lemmas \ref{l.linear}, \ref{l.linearII} and (\ref{e.global.9}), it follows that
$$\|\psi_1(T)\|_{H^{-1/2}}\leq \|\psi_{10}\|_{H^{-1/2}} + C\cdot T^{\tfrac{1}{2}-3\varepsilon}\|B_0\|_{L^2}^2$$
and
$$\|\psi_2(T)\|_{H^{-1/2}}\leq \|\psi_{20}\|_{H^{-1/2}} + C\cdot T^{\tfrac{1}{2}-3\varepsilon}\|B_0\|_{L^2}^2$$
since
$$\psi_1(t)=W_+(\psi_{10})(t)+W_+*_R(\partial_x(|B|^2)) \textrm{ and }
\psi_2(t)=W_-(\psi_{20})(t)-W_-*_R(\partial_x(|B|^2)).$$
Thus, we can perform $m$ iterations of this scheme before the quantity $\|\psi_1(t)\|_{H^{-1/2}}$ doubles, where
\begin{equation}\label{e.global.11}
m\sim \frac{\|\psi_{10}\|_{H^{-1/2}}}{T^{\tfrac{1}{2}-3\varepsilon}\|B_0\|_{L^2}^2}.
\end{equation}
On the other hand, during these $m$ iterations, either $\|\psi_2(t)\|_{H^{-1/2}}\lesssim \|B_0\|_{L^2}^2$
for all $t\in [0,mT]$ or $\|\psi_2(jT)\|_{H^{-1/2}}\gg \|B_0\|_{L^2}^2$ at some stage $1<j<m$:
\begin{itemize}
\item in the first situation, it is clear that this iteration scheme allows us to produce a solution of the
Zakharov-Rubenchik system on the time interval $[0,mT]$ with
$$mT\sim\|B_0\|_{L^2}^{-2}$$
by (\ref{e.global.10}) and (\ref{e.global.11});
\item in the second situation, we observe that $(B(jT),\psi_1(jT),\psi_2(jT))$ fits the assumptions of the
previous case a).
\end{itemize}
Therefore, we are able to conclude (in any situation) the global well-posedness result in the case b) and the estimate
$$\|\psi_1(t)\|_{H^{-1/2}}\leq \exp(Ct\|B_0\|_{L^2}^2)\max\{\|\psi_{10}\|_{H^{-1/2}}, \|B_0\|_{L^2}^2\}$$
and
$$\|\psi_2(t)\|_{H^{-1/2}}\leq \exp(Ct\|B_0\|_{L^2}^2)\max\{\|\psi_{20}\|_{H^{-1/2}}, \|B_0\|_{L^2}^2\}.$$
\end{itemize}
This shows that the global well-posedness result of Theorem~\ref{t.global} for
initial data $(B_0,\psi_{10},\psi_{20})$ in $L^2\times H^{-1/2}\times
H^{-1/2}$ holds and, moreover, we have the estimates
\begin{equation}\label{e.global.12}
\|\psi_1(t)\|_{H^{-1/2}}\leq
\exp(Ct\|B_0\|_{L^2}^2)\max\{\|\psi_{10}\|_{H^{-1/2}},
\|B_0\|_{L^2}^2\}
\end{equation}
and
\begin{equation}\label{e.global.13}
\|\psi_2(t)\|_{H^{-1/2}}\leq \exp(Ct\|B_0\|_{L^2}^2)\max\{\|\psi_{20}\|_{H^{-1/2}}, \|B_0\|_{L^2}^2\}.
\end{equation}

\bigskip

\bigskip

Once the case $k=0$ and $l=-1/2$ of Theorem~\ref{t.global} is
proved, we can deal with \eqref{zrnew} for initial data
$(B_0,\psi_{10},\psi_{20})\in H^k\times H^{l}\times H^{l}$ (where
$0\leq k=l+\tfrac{1}{2}$) as follows: estimating (\ref{e.global.1})
in $X^{k,b}\times W_+^{l,b_1}\times W_-^{l,b_1}$ (for $b=3/8+$,
$b_1=1/4+$), we obtain
\begin{equation}\label{e.global.14}
\begin{split}
&\|\Phi_T(B,\psi_1,\psi_2)\|_{X^{k,b}}\lesssim
T^{\tfrac{1}{8}-}\|B_0\|_{H^k} +
T^{\tfrac{1}{4}-}(\|\psi_1\|_{W_+^{-1/2,b_1}}+
\|\psi_2\|_{W_-^{-1/2,b_1}}+\|B\|_{X^{0,b}}^2)\|B\|_{X^{k,b}} \\
&+T^{\tfrac{1}{4}-}\|B\|_{X^{0,b}}(\|\psi_1\|_{W_+^{l,b_1}}+
\|\psi_2\|_{W_-^{l,b_1}})\\
&\|\Psi_T^+(B)\|_{W_+^{l,b_1}}\lesssim
T^{\tfrac{1}{4}-}\|\psi_{10}\|_{H^{l}} +
T^{\tfrac{1}{2}-}\|B\|_{X^{0,b}}\|B\|_{X^{k,b}} \\
&\|\Psi_T^-(B)\|_{W_-^{-1/2,b_1}}\lesssim
T^{\tfrac{1}{4}-}\|\psi_{20}\|_{H^{l}} +
T^{\tfrac{1}{2}-}\|B\|_{X^{0,b}}\|B\|_{X^{k,b}}
\end{split}
\end{equation}
by the Lemmas~\ref{l.linear},~\ref{l.linearII},~\ref{l.multilinear-2}. On
the other hand, since $0\leq k=l+\tfrac{1}{2}$, we know that
$(B_0,\psi_{10},\psi_{20})\in L^2\times H^{-1/2}\times H^{-1/2}$ and
the corresponding solution $(B,\psi_1,\psi_2)$ satisfy the \emph{a
priori} estimates (\ref{e.global.12}) and (\ref{e.global.13})
(besides the upper bound $\|B\|_{X^{0,b}}\lesssim \|B_0\|_{L^2}$) on
a time interval $[0,T]$ for some $T=T(\|B_0\|_{L^2}^{-2})$. Thus,
using these facts into (\ref{e.global.14}), we get uniform estimates
for $\|B\|_{X^{k,b}}$, $\|\psi_1\|_{W_+^{l,b_1}}$ and
$\|\psi_2\|_{W_-^{l,b_1}}$ on a time interval of length
$T=T(\|B_0\|_{L^2}^{-2})$. This completes the proof.
\end{proof}

\section{Polynomial growth of higher Sobolev norms}

The goal of this section is the study of certain upper bounds for the
higher Sobolev norms of the solutions $(B,\psi_1,\psi_2)$ of the
Zakharov-Rubenchik system~\eqref{zrnew}. The strategy adopted here follows the lines of
Colliander and Staffilani~\cite{CS}. In particular, we construct global-in-time solutions
of~\eqref{zrnew} for initial data $(B_0,\psi_{10},\psi_{20})\in H^1\times L^2\times L^2$ and
we show that the Sobolev norm $\|B\|_{H^s}$ of the Schr\"odinger part $B$ of Zakharov-Rubenchik system
satisfy certain \emph{polynomial} upper bounds for all $s\gg 1$. More precisely, we prove the following
result:

\begin{theorem}\label{t.polynomial}The Zakharov-Rubenchik~\eqref{zrnew} is globally well-posed for initial data $$(B_0,\psi_{10},\psi_{20})\in H^1\times L^2\times L^2.$$
Furthermore, if $B_0,\psi_{10},\psi_{20}\in \mathcal{S}$, where
$\mathcal{S}$ is the Schwartz class, then this global solution
$(B,\psi_1,\psi_2)$ satisfy
\begin{equation}\label{e.polynomial}
\|B(t)\|_{H^s}\lesssim 1+|t|^{(s-1)+}.
\end{equation}
\end{theorem}

\begin{remark}Since we can write
\begin{equation*}
\begin{cases}
\psi_1(t)=W_+(t)\psi_{10}+W_+*_R\partial_x(|B|^2)(t)\\
\psi_2(t)=W_-(t)\psi_{20}-W_-*_R\partial_x(|B|^2)(t),
\end{cases}
\end{equation*}
it is not hard to infer regularity bounds (in $H^{s-1}$) for $\psi_1$ and $\psi_2$ from the estimate
\eqref{e.polynomial}.
\end{remark}

In the sequel, we subdivide this section into two parts: the first subsection is dedicated to the global
well-posedness statement of Theorem~\ref{t.polynomial} and the second subsection contains the proof of
the estimate~\eqref{e.polynomial}.

\subsection{Global well-posedness in $H^1(\R)\times L^2(\R)\times L^2(\R)$}

Using the conservation laws $I_1(t)$, $I_2(t)$, $I_3(t)$ and $I_4(t)$ in \eqref{ene1}-\eqref{ene4}
the following \emph{a priori} estimate can be established.

\begin{lemma}\label{l.apriori} Solutions of the IVP associated to \eqref{zrnew} satisfy
\begin{equation*}
\begin{split}
\|(B(t),\psi_1(t),\psi_2(t))\|_{H^1\times L^2\times L^2}
\lesssim \|(B_0,\psi_{10},\psi_{20})\|_{H^1\times L^2\times L^2}^2+ I_1(0)^3.
\end{split}
\end{equation*}
\end{lemma}

\begin{proof} It follows using the conserved quantities (see Lemma 3.3 of~\cite{O}) and \eqref{changevar}.
\end{proof}

On the other hand, the local well-posedness result of Theorem~\ref{t.local} ensures the existence
of a local-in-time solution $(B,\psi_1,\psi_2)$ for initial data
$(B_0,\psi_{10},\psi_{20})\in H^1\times L^2\times L^2$
in a time interval $[0,T]$, where
\begin{equation}\label{e.localpoly}
T > c\|(B_0,\psi_{10},\psi_{20})\|_{H^1\times L^2\times L^2}^{-\alpha}
\end{equation}
for some $\alpha>0$. Combining this fact with Lemma~\ref{l.apriori}, we obtain the following result:

\begin{proposition}\label{p.globalpoly}The Zakharov-Rubenchik system~\eqref{zrnew} is globally well-posed
in $H^1\times L^2\times L^2$.
\end{proposition}

\subsection{Polynomial upper bounds for $\|B(t)\|_{H^s}$} In view of the Proposition \ref{p.globalpoly}
above, the proof
of Theorem~\ref{t.polynomial} is complete once we show the following result:

\begin{proposition}\label{p.polybound1}For initial data $B_0,\psi_{10},\psi_{20}\in\mathcal{S}$, the global
solution of \eqref{zrnew} satisfy
$$
\|B(t)\|_{H^s}\lesssim 1+|t|^{(s-1)+}.
$$
\end{proposition}

\begin{proof}We begin by noticing that it suffices to bound $\|B(t)\|_{H^s}$ for $t\in [0,T]$ (with T
defined by \eqref{e.localpoly}). Since $\|B(t)\|_{L^2}=\|B_0\|_{L^2}$ for all $t$, it remains only to
compute $\|A^sB(t)\|_{L^2}$, where $A^s:=\sqrt{-\Delta}$. To avoid technical issues involving fractional
derivatives, we assume $s=2m$ a large even integer ($m\gg 1$). Denoting by $\langle.,.\rangle$ the usual
$L^2$ inner product, we obtain
\begin{equation*}
\begin{split}
\|A^sB(t)\|_{L^2}&=\|A^sB_0\|_{L^2}+\int_0^t\frac{d}{d\mu}\langle A^sB(\mu),A^sB(\mu)\rangle d\mu\\
&= \|A^sB_0\|_{L^2}+\Re\int_0^t\frac{d}{d\mu}\langle A^sB(\mu),A^sB(\mu)\rangle d\mu.
\end{split}
\end{equation*}
This reduces matters to the computation of the term
\begin{equation*}
I:=\Re\int_0^t\frac{d}{d\mu}\langle A^sB(\mu),A^sB(\mu)\rangle d\mu =
2\Re\int_0^t\langle A^s\dot{B}(\mu),A^sB(\mu)\rangle d\mu.
\end{equation*}
Because the Schr\"odinger part $B$ of the Zakharov-Rubenchik system~\eqref{zr1d}
satisfies\footnote{Here we are omitting the irrelevant constants since our task is to show a polynomial
upper bound for $\|B(t)\|_{H^s}$.}
\begin{equation}\label{e.Duhamel-B}
\begin{split}
\partial_t B &= i\Delta B-i\psi_1 B + i\psi_2 B-i |B|^2 B \\
&=i\Delta B - i \left(W_+(\psi_{10})-W_-(\psi_{20})\right)B
-i \left(W_+*_R\partial_x(|B|^2)-W_-*_R\partial_x(|B|^2)\right)B -i |B|^2B,
\end{split}
\end{equation}
we conclude that
\begin{equation*}
\begin{split}
I&=-2\Im\int_0^t\langle A^s \Delta B(\mu), A^s B(\mu)\rangle d\mu +
2\Im\int_0^t \langle A^s\left(W_{\pm}(\psi_{10}, \psi_{20})(\mu)B(\mu)\right), A^s B(\mu)\rangle d\mu \\
&+ 2 \Im\int_0^t \langle A^s\left(W_{\pm}*_R \partial_x(|B|^2)(\mu)B(\mu)\right), A^s B(\mu)\rangle d\mu +
2\Im\int_0^t \langle A^s\left(|B(\mu)|^2 B(\mu)\right),A^s B(\mu)\rangle d\mu \\
&:=I_1+I_2+I_3+I_4,
\end{split}
\end{equation*}
where
$$
W_{\pm}(\psi_{10},\psi_{20})(t):= W_+(\psi_{10})(t)-W_-(\psi_{20})(t)
$$
and
$$
W_{\pm}*_R \partial_x(|B|^2)(t):= W_+*_R\partial_x(|B|^2)(t)-W_-*_R\partial_x(|B|^2)(t).
$$
Consider the term $I_1$. Since $-\Delta=A^2$, we get from integration by parts that the integrand of $I_1$ is
a real number, so that $I_1=0$. On the other hand, the term $I_2$ involves the expression
$A^s(W_{\pm}(\psi_{10},\psi_{20})\cdot B)$. This expression can be expanded
using the Leibnitz rule. The term with the highest derivative on $B$ is $W_{\pm}(\psi_{10},\psi_{20})A^s B$,
but it does not contribute for the computation of $I_2$ since $W_{\pm}(\psi_{10},\psi_{20})$ is a real-valued
function so that the corresponding integrand is a real number. Thus, $I_2$ becomes a sum of terms of the form:
\begin{equation*}
C\Im\int_0^t\langle A^{s_1}W_{\pm}(\psi_{10},\psi_{20})(\mu)\cdot A^{s_2} B(\mu), A^s B(\mu)\rangle d\mu
\end{equation*}
where $s=s_1+s_2$, $s_1,s_2\in\mathbb{N}$, $1\leq s_1\leq s$ and $0\leq s_2\leq s-1$. Multiplying by a smooth
localized-in-time cutoff function $\kappa_T$ supported on the interval $[0,T]$ and using H\"older inequality,
we can estimate each of these terms by
\begin{equation*}
\|A^{s_1}W_{\pm}(\psi_{10},\psi_{20})\|_{L^2_{xt}}\|A^{s_2}B\|_{L^4_{xt}}\|A^s B\|_{L^4_{xt}}.
\end{equation*}
Applying the Bourgain-Strichartz estimate $\|f\|_{L^4_{xt}}\lesssim \|f\|_{X^{0,3/8+}}$ and the local
well-posedness result $\|B\|_{X^{m,b}}\lesssim \|B_0\|_{H^m}$ (for $b=1/2+$ and $m\geq 0$), we obtain
that $I_2$ can be estimated by a sum of terms of the form
\begin{equation*}
(\|\psi_{10}\|_{H^{s_1}}+\|\psi_{20}\|_{H^{s_1}})\|B_0\|_{H^{s_2}}\|B_0\|_{H^s}.
\end{equation*}
Interpolating $\|B_0\|_{H^{s_2}}$ between $\|B_0\|_{H^1}$ and $\|B_0\|_{H^s}$, we get the following inequality:
\begin{equation}\label{e.I2}
|I_2|\lesssim \sum\limits_{s_2=1}^{s-1}\|B_0\|_{H^s}^{1+\frac{s_2-1}{s-1}}
\lesssim \|B_0\|_{H^s}^{2-\frac{1}{s-1}}.
\end{equation}
Similarly, $I_4$ is a sum of terms of the form:
\begin{equation*}
I_4(s_0,s_1,s_2):=C\Im\int_0^t\langle A^{s_0}B(\mu)\overline{A^{s_1}B(\mu)}A^{s_2}B(\mu), A^s B(\mu)\rangle d\mu
\end{equation*}
where $s_0+s_1+s_2=s$, $s_0,s_1,s_2\in\mathbb{N}$. Since the integrand of $I_4(s,0,0)$ and $I_4(0,0,s)$,
we have $I_4(s,0,0)=I_4(0,0,s)=0$. On the other hand, the bilinear estimates
provided by Lemmas 3.3, 3.4 and 3.7 of~\cite{S} give the estimate
\begin{equation*}
|I_4(0,s,0)|\lesssim \|B\|_{X^{1,1/2+}}^2 \|B\|_{X^{s-3/4+,1/2+}}^2.
\end{equation*}
Hence, from the local well-posedness theory, we deduce that
$$
|I_4(0,s,0)|\lesssim \|B_0\|_{H^1}^2\|B_0\|_{H^{s-3/4+}}^2.
$$
Next, we estimate $I_4(s_0,s_1,s_2)$
for $1\leq s_0,s_1,s_2\leq s-1$, $s_0+s_1+s_2=s$. Denoting by $j_2=\max\{s_0,s_1,s_2\}$,
$j_0=\min\{s_0,s_1,s_2\}$ and $j_1\in\{s_0,s_1,s_2\}-\{j_0,j_2\}$ and using again the Lemmas 3.3
and 3.4 of~\cite{S} (plus the local well-posedness theory), we obtain
\begin{equation*}
\begin{split}
|I_4(s_0,s_1,s_2)|&\lesssim\|B\|_{X^{j_0+1,1/2+}}\|B\|_{X^{j_2+1,1/2+}}\|B\|_{X^{j_2-3/4+,1/2+}}\|B\|_{X^{s-3/4+,1/2+}} \\
&\lesssim \|B_0\|_{H^{j_0+1}}\|B_0\|_{H^{j_1+1}}\|B_0\|_{H^{j_2-3/4+}}\|B_0\|_{H^{s-3/4+}}.
\end{split}
\end{equation*}
Therefore, summing up the bounds for the terms $I_4(s_0,s_1,s_2)$ and interpolating
$\|B_0\|_{H^{\widetilde{s}}}$ between $H^1$ and $H^s$, we see that
\begin{equation}\label{e.I4}
|I_4|\lesssim \|B_0\|_{H^s}^{(2s-2-3/2+)/(s-1)}.
\end{equation}
Finally, it remains to deal with $I_3$. Using the Leibnitz rule, we can write as a sum of terms of the form
\begin{equation*}
I_3(s_0,s_1,s_2)=C\int_0^t\langle W_{\pm}*_R(\partial_x(A^{s_0}B(\mu)A^{s_1}\overline{B}(\mu)))\cdot A^{s_2}B(\mu), A^sB(\mu)\rangle d\mu
\end{equation*}
where $s_0+s_1+s_2=s$. Note that the integrand of $I_3(0,0,s)$ is real so that $I_3(0,0,s)=0$. Now we treat
the terms $I_3(s_0,s_1,s_2)$ with $s_0+s_1+s_2=s$ and $s_2\leq s-1$. By Cauchy-Schwarz inequality we get
\begin{equation*}
|I_3(s_0,s_1,s_2)|\lesssim \left\|W_{\pm}*_R(\partial_x(A^{s_0}B(\mu)A^{s_1}\overline{B}(\mu)))\right\|_{L_{xt}^2} \left\|A^{s_2}B(\mu) A^s\overline{B}(\mu)\right\|_{L_{xt}^2}
\end{equation*}
Using the bilinear estimate of Lemma~\ref{l.multilinear} (with $l=0$, $k=1/4$), the identity $L_{xt}^2=W_{\pm}^{0,0}$,
Bourgain's refinement of Strichartz estimate~\cite{B} and the local well-posedness theory, it follows that
\begin{equation*}
\begin{split}
|I_3(s_0,s_1,s_2)|&\lesssim \|B\|_{X^{s_0+1/4,1/2+}} \|B\|_{X^{s_1+1/4,1/2+}}\|B\|_{X^{s-1/2+,1/2+}}^2\\
&\lesssim \|B_0\|_{H^{s_0+1/4}}\|B_0\|_{H^{s_1+1/4}}\|B_0\|_{H^{s_3-1/2+}}\|B_0\|_{H^{s-1/2}}.
\end{split}
\end{equation*}
By interpolation, we infer that
\begin{equation}\label{e.I3}
|I_3(s_0,s_1,s_2)|\lesssim \|B_0\|_{H^s}^{(2s-2-5/2+)/(s-1)}.
\end{equation}
Combining the estimates~\eqref{e.I2},~\eqref{e.I4} and~\eqref{e.I2} with the fact $I_1=0$, we obtain that
\begin{equation}
\|B(t)\|_{H^s}\leq \|B_0\|_{H^s} + C\|B_0\|_{H^s}^{1-\tfrac{1}{s-1}+}.
\end{equation}
This ends the proof of Proposition~\ref{p.polybound1}.
\end{proof}

\section{Ill-posedness of the Zakharov-Rubenchik}

In this section we establish Theorems \ref{t.ill-inflation}-\ref{t.ill-decoherence}. They are concerned
with the lack of smoothness of the Zakharov-Rubenchik evolution for initial data
$(B_0,\psi_{10},\psi_{20})\in H^k\times H^l\times H^s$
when $(k,l,s)$ stays outside the local well-posedness region established in Theorem~\ref{t.local}.

\subsection{Proof of Theorem~\ref{t.ill-inflation}} In the sequel, we follow closely the ``norm-inflation'' scheme
for the Zakharov system discussed in \cite{H}. In particular, since the evolutions of the Zakharov system and
the Zakharov-Rubenchik system are somewhat similar\footnote{In the sense that the linear groups involved in the
analysis of both systems (namely Schr\"odinger and transport) are equal.}, we will just point out certain
modifications of the arguments of section 3 of~\cite{H} in order to get the result of Theorem~\ref{t.ill-inflation}.

Let $0<k<1$, $l>2k-\tfrac{1}{2}$. Firstly, we suppose also that $l$ is \emph{near} $2k-\tfrac{1}{2}$, i.e.,
\begin{itemize}
\item If $0<k\leq \tfrac{1}{4}$, $l$ is restricted to the region $2k-\tfrac{1}{2}<l\leq 4k-\tfrac{1}{2}$;
\item If $\tfrac{1}{4}<k<1$, $l$ is restricted to the region $2k-\tfrac{1}{2}<l\leq \tfrac{4k}{3}+\tfrac{1}{6}$.
\end{itemize}
Under these assumptions, we fix
$$\sigma:=
\begin{cases}
k, \quad \quad \quad \textrm{ if } 0<k\leq\tfrac{1}{4} \textrm{ and } 2k-\tfrac{1}{2}<l\leq 4k-\tfrac{1}{2} \\
\tfrac{k}{3}+\tfrac{1}{6}, \quad \, \textrm{ if } \tfrac{1}{4}<k<1 \textrm{ and }
2k-\tfrac{1}{2}<l\leq \tfrac{4k}{3}+\tfrac{1}{6}.\\
\end{cases}$$
Next, we put $f_N=f_{N,A}+f_{N,B}$ where
$$\widehat{f_{N,A}}(\xi):=N^{\tfrac{1}{2}-k}\chi_{[-N-\tfrac{1}{N},-N]}(\xi) \quad \textrm{and} \quad
\widehat{f_{N,B}}(\xi):=N^{\tfrac{1}{2}-k}\chi_{[N+1,N+1+\tfrac{1}{N}]}(\xi).$$
Recall that the solution $(B_N,\psi_{1,N},\psi_{2,N})$ of the Zakharov-Rubenchik system with
initial data $(f_N,0,0)$ satisfies
$$
\begin{cases}
B_N(t)=\kappa(t)U(t)f_N + i \kappa_T(t)U*_R\{|B_N|^2 + (W_+-W_-)*_R\partial_x(|B_N|^2)\}\cdot B_N(t),\\
\psi_{1,N}(t)=\kappa_T(t)W_+*_R\partial_x(|B_N|^2)(t),\\
\psi_{2,N}(t)=\kappa_T(t)W_-*_R\partial_x(|B_N|^2)(t).
\end{cases}
$$
During the section 3 of~\cite{H}, Holmer proved
$$\|W_+*_R\partial_x(|Uf_N|^2)(t)\|_{H^l}\sim t N^{l-(2k-\tfrac{1}{2})}$$
for $N\gtrsim t^{-1}$ (see estimate (3.7) of~\cite{H}). On the other hand, the linear and multilinear
estimates of Lemma~\ref{l.multilinear} can be used in the same manner as Holmer~\cite{H} (p.10, 11) to
give us that
\begin{eqnarray*}
\|B_N-\kappa(t)U(t)f_N\|_{X^{k+\sigma,b_1}}&\lesssim&
\|W_+*_R\partial_x(|B_N|^2) \cdot B_N\|_{X^{k+\sigma,b_1}}+\\
&+& \|W_-*_R\partial_x(|B_N|^2) \cdot B_N\|_{X^{k+\sigma,b_1}} + \| |B_N|^2 B_N\|_{X^{k+\sigma,b_1}} \\
&\lesssim& \|f_N\|_{H^{k'}}^2\|f_N\|_{H^k}\sim N^{2(k'-k)+\sigma},
\end{eqnarray*}
where $k'=0$, $b_1=\tfrac{3}{4}-\tfrac{k+\sigma}{2}$ if $0<k+\sigma<\tfrac{1}{2}$ and
$k'=\tfrac{k+\sigma}{2} - \tfrac{1}{4}$, $b_1=\tfrac{1}{2}$ if $\tfrac{1}{2}\leq k+\sigma<\tfrac{5}{2}$.
Combining these estimates and reasoning as Holmer~\cite{H} (p.11, 12) lead us to
$$\|\psi_{1,N}(t)\|_{H^l}\gtrsim t\cdot N^{l-(2k-\tfrac{1}{2})}$$
for $N\gtrsim t^{-1}$. Therefore, we showed the first part of Theorem~\ref{t.ill-inflation} under the
assumptions $0<k<1$, $l>2k-\tfrac{1}{2}$ and $l$ near $2k-\tfrac{1}{2}$ (in the sense above). Finally,
the general case, i.e., either $0<k<1$, $l>2k-\tfrac{1}{2}$ or $k<0$, $l>-1/2$ can be reduced to this
previous particular case by decreasing $l$ and increasing $k$ appropriately (since
$\|F\|_{H^{q'}}\leq \|F\|_{H^q}$ whenever $q'\leq q$).

Analogously, when dealing with the second part of Theorem~\ref{t.ill-inflation}, it suffices to replace
$l$ by $s$ and $f_N$ by $g_N=g_{N,A}+g_{N,B}$ where
$$\widehat{g_{N,A}}(\xi):=N^{\tfrac{1}{2}-k}\chi_{[-N-\tfrac{1}{N},-N]}(\xi) \quad \textrm{and} \quad
\widehat{g_{N,B}}(\xi):=N^{\tfrac{1}{2}-k}\chi_{[N-1,N-1+\tfrac{1}{N}]}(\xi)$$
in the previous considerations so that the same argument applies. The details are left to the reader.
This completes the proof of Theorem~\ref{t.ill-inflation}.

\subsection{Proof of Theorem~\ref{t.ill-C2}}For sake of simplicity, we assume that $k\in\mathbb{R}$
and $l<-1/2$ (since the case $s<-1/2$ is similar). We take $N$ a large integer and we put
$\widehat{B_0}(\xi):=N^{\tfrac{1}{2}-k}\chi_{[0,\tfrac{1}{N}]}(\xi)$,
$\widehat{\psi_{10}}(\xi):=N^{\tfrac{1}{2}-l} \chi_{[-\tfrac{1}{N},\tfrac{1}{N}]}(\xi)$ and
$\psi_{20}\equiv 0$. Let $\gamma\in\mathbb{R}$ be a parameter and denote by
$F(B_0,\psi_{10},\psi_{20})=(B(t),\psi_1(t),\psi_2(t))$ the data-to-solution map of the
Zakharov-Rubenchik system. Assume that $F$ is $C^2$ at the origin $(0,0,0)$ and consider the
path $G(\gamma)=(\gamma B_0,\gamma\psi_{10},\gamma\psi_{20}) = (\gamma B_0,\gamma\psi_{10},0)$.

Note that the solution $(B,\psi_1,\psi_2)=F\circ G(\gamma)$ of~\eqref{zrnew} verifies
$$
\begin{cases}
B(t)=U(t)(\gamma B_0) + i \int_0^t U(t-t')\{|B|^2 + \psi_1+\psi_2\}\cdot B(t') dt',\\
\psi_{1}(t)=W_+(t)(\gamma\psi_{10})+\int_0^t W_+(t-t')\partial_x(|B|^2)(t')dt',\\
\psi_{2}(t)=-\int_0^t W_-(t-t')\partial_x(|B|^2)(t')dt'.
\end{cases}
$$

Thus, we get
\begin{eqnarray*}
\partial_\gamma B(t)&=&U(t)B_0 + i \int_0^t U(t-t')\{\partial_\gamma B(\psi_1+\psi_2+|B|^2) + \\
&+&B(\partial_\gamma\psi_1 +\partial_\gamma\psi_2 + \partial_\gamma B\cdot \overline{B} + B\overline{\partial_\gamma B})\}(t') dt'
\end{eqnarray*}
and
$$
\begin{cases}
\partial_\gamma\psi_{1}(t)=W_+(t)\psi_{10}+\int_0^t W_+(t-t')\partial_x(\partial_\gamma B\cdot\overline{B} + B\cdot \overline{\partial_\gamma B})(t')dt',\\
\partial_\gamma\psi_{2}(t)=-\int_0^t W_-(t-t')\partial_x(\partial_\gamma B\cdot\overline{B} + B\cdot \overline{\partial_\gamma B})(t')dt'.
\end{cases}
$$
Using that $F\circ G(0)=(0,0,0)$, it follows
$$
\begin{cases}
\partial_\gamma B|_{\gamma=0}(t)=U(t)B_0\\
\partial_\gamma\psi_{1}(t)=W_+(t)\psi_{10},\\
\partial_\gamma\psi_{2}(t)=0.
\end{cases}
$$
Taking derivative with respect to $\gamma$ leads us to
\begin{eqnarray*}
\partial^2_\gamma B(t)&=& i \int_0^t U(t-t')\{\partial^2_\gamma B(\psi_1+\psi_2+|B|^2)
+ 2\partial_\gamma B(\partial_\gamma\psi_1
+\partial_\gamma\psi_2+\partial_\gamma B\cdot\overline{B} + B \overline{\partial_\gamma B})\\
&+&B(\partial^2_\gamma\psi_1 +\partial^2_\gamma\psi_2 + \partial^2_\gamma B\cdot \overline{B} +
2 |\partial_\gamma B|^2 + B\overline{\partial^2_\gamma B})\}(t') dt'
\end{eqnarray*}
and
$$
\begin{cases}
\partial^2_\gamma\psi_{1}(t)=\int_0^t W_+(t-t')\partial_x(\partial^2_\gamma B\cdot\overline{B} +
2|\partial_\gamma B|^2 + B\cdot\overline{\partial^2_\gamma B})(t')dt',\\
\partial^2_\gamma\psi_{2}(t)=-\int_0^t W_-(t-t')\partial_x(\partial^2_\gamma B\cdot\overline{B} +
2|\partial_\gamma B|^2 + B\cdot\overline{\partial^2_\gamma B})(t')dt'.
\end{cases}
$$
Hence,
$$\partial^2_\gamma B|_{\gamma=0}(t) = 2i\int_0^t U(t-t') \{U(t')B_0 \cdot W_+(t')\psi_{10}\} dt'.$$
If $F$ is $C^2$ at the origin, we can use that $D^2F(0,0,0)$ is a bounded bilinear operator to conclude
that the following bilinear estimate holds:
$$
\|\int_0^t U(t-t') \{U(t')B_0 \cdot W_+(t')\psi_{10}\} dt'\|_{H^k}
\lesssim\|B_0\|_{H^k}\|\psi_{10}\|_{H^l}.
$$
Since $\|B_0\|_{H^k}=\|\psi_{10}\|_{H^l}=1$, we get
\begin{equation}\label{e.C2-1}
\|\int_0^t U(t-t') \{U(t')B_0 \cdot W_+(t')\psi_{10}\} dt'\|_{H^k}\lesssim 1.
\end{equation}
On the other hand, we know that $L(x,t):=\int_0^t U(t-t') \{U(t')B_0 \cdot W_+(t')\psi_{10}\} dt'$ satisfies
$$
\widehat{L}(\xi,t) = \int_0^t e^{-i(t-t')\xi^2} \int e^{-it'\xi_1^2}\widehat{B_0}(\xi_1) e^{-it'(\xi-\xi_1)}\widehat{\psi_{10}}(\xi-\xi_1) d\xi_1dt'.
$$
Thus,
$$\widehat{L}(\xi,t) = e^{-it\xi^2} \int \widehat{B_0}(\xi_1)\widehat{\psi_{10}}(\xi-\xi_1) \frac{e^{it(\xi-\xi_1)(\xi+\xi_1-1) -1}}{ i(\xi-\xi_1)(\xi+\xi_1-1)} d\xi_1.$$
Because the support of $\widehat{B_0}(\xi_1)\widehat{\psi_{10}}(\xi-\xi_1)$ is contained
in the region $\xi_1\in[0,\tfrac{1}{N}]$ and $\xi\in[-\tfrac{1}{N},\tfrac{2}{N}]$, we can use the expansion
of $e^z$ to obtain
$$
\widehat{L}(\xi,t) = e^{-it\xi^2} (t+O(t^2)) N^{-l-k}.
$$
In particular,
\begin{equation}\label{e.C2-2}
\|\int_0^t U(t-t') \{U(t')B_0 \cdot W_+(t')\psi_{10}\} dt'\|_{H^k}\gtrsim t\cdot N^{-l-\tfrac{1}{2}}.
\end{equation}
Combining the estimates~\eqref{e.C2-1} and~\eqref{e.C2-2}, we get a contradiction for $N$ sufficiently
large (since $l<-1/2$). This ends the proof of Theorem~\ref{t.ill-C2}.

\subsection{Proof of Theorem~\ref{t.ill-decoherence}}
The idea is to find a suitable class of initial data so that the Schr\"odinger variables of the
Zakharov-Rubenchik evolution eventually exhibit completely different phases.

For the sake of simplicity, we will normalize the constants depending on $\omega,k,\beta,\nu,q,\theta$ from
the system~\eqref{zrnew} so that the ZR system becomes
\begin{equation}\label{e.zr'}
\begin{cases}
i\partial_t B + \partial_x^2 B = \psi_+ B + \psi_- B + |B|^2 B, \\
\partial_{t}\psi_+ + \partial_{x}\psi_+ = \partial_{x}(|B|^2),\\
\partial_t\psi_- -\partial_x\psi_- = \partial_x(|B|^2).
\end{cases}
\end{equation}
Also, up to exchange the roles of $\psi_+$ and $\psi_-$ in the arguments below, it suffices to prove the
Theorem~\ref{t.ill-decoherence} in the case $k=0$ and $l<-3/2$.

Next, we fix four parameters $L\gg 1$, $0<\mu\ll 1$, $0<c<1$ and $0<\Theta\ll 1$. Consider the following modified ZR system:
\begin{equation}\label{e.mzr'}
\begin{cases}
i\partial_t \widetilde{B} + \mu^2\partial_x^2 \widetilde{B} = \widetilde{\psi}_{+0}\left(x-\tfrac{\mu(1-c)}{L}t\right) \widetilde{B} + \widetilde{\psi}_{-0}\left(x+\tfrac{\mu(1+c)}{L}t\right) \widetilde{B}
+ (\widetilde{\psi}_+ +\widetilde{\psi}_-) \widetilde{B} + \Theta^2 |\widetilde{B}|^2 \widetilde{B}, \\
\tfrac{L}{\mu(1-c)}\partial_{t}\widetilde{\psi}_+ + \partial_{x}\widetilde{\psi}_+ = \tfrac{\Theta^2}{1-c}\partial_{x}(|\widetilde{B}|^2),\\
\tfrac{L}{\mu(1+c)}\partial_t\widetilde{\psi}_- - \partial_x\widetilde{\psi}_- = \tfrac{\Theta^2}{1+c}\partial_x(|\widetilde{B}|^2),
\end{cases}
\end{equation}
with initial conditions
\begin{equation*}
\begin{cases}
\widetilde{B}(0,x) = \widetilde{B}_0(x), \\
\widetilde{\psi}_+(0,x)=0,
\\
\widetilde{\psi}_-(0,x) = 0.
\end{cases}
\end{equation*}

\begin{proposition}\label{p.step2}If $k\geq 1$ and
\begin{equation}\label{e.c1}
T_0\lesssim |\log\mu|, \quad L\gtrsim \mu^{-5}, \quad \Theta^2 \sim \mu L^{-1},
\end{equation}
then it holds
\begin{equation}\label{e.s-d-1}
\|\widetilde{B}\|_{L_{[0,T_0]}^{\infty}H_x^k}\lesssim\mu^{-1/2},
\end{equation}
\begin{equation}\label{e.s-d-2}
\|\widetilde{\psi}_{\pm}\|_{L_{[0,T_0]}^{\infty}H_x^{k-1}}\lesssim\frac{\Theta^2}{L}.
\end{equation}
\end{proposition}
\begin{proof}Since $\psi_{\pm} = \tfrac{\Theta^2}{(1\mp c)}\int_0^{t\alpha_\pm}\partial_x(|\widetilde{B}|^2)(x-t', t-\tfrac{t'}{\alpha_\pm}) dt'$,
where $\alpha_\pm = L/\mu(1\mp c)$, it follows that
\begin{equation}\label{e.s-d-3}
\|\widetilde{\psi}_{\pm}\|_{L_{[0,T_0]}^{\infty}H_x^{k-1}}\leq \frac{\Theta^2\mu T}{L}
\|\widetilde{B}\|_{L_{[0,T_0]}^{\infty}H_x^k}^2.
\end{equation}
On the other hand, using the energy method for $\widetilde{B}$, we get
\begin{eqnarray*}
&\|\partial_x^k\widetilde{B}(T_0)\|_{L^2_x}^2 - \|\partial_x^k\widetilde{B}(0)\|_{L^2_x}^2 \\
&=-\Re\left(i\int_0^{T_0}\int\partial_x^k(\psi_{\pm0}(...)\widetilde{B})\overline{\partial_x^k\widetilde{B}} dx dt\right)
- \Re\left(i\int_0^{T_0}\int\partial_x^k(\psi_{\pm}\widetilde{B})\overline{\partial_x^k\widetilde{B}} dx dt\right)  \\
&-\Theta^2 \Re\left(i\int_0^{T_0}\int\partial_x^k(|\widetilde{B}|^2\widetilde{B})\overline{\partial_x^k\widetilde{B}} dx dt\right) \\
&= (I) + (II) + (III).
\end{eqnarray*}
Note that,
$$
|(I)|\leq \|\psi_{\pm0}\|_{H^k}\int_0^{T_0}\|\widetilde{B}(t)\|_{H^k}^2 dt.
$$
Now, following the calculations of Holmer~\cite{H} (p.16), it is not hard to see that
$$
|(II)|\lesssim \left(\frac{\Theta^2\mu T^2}{L}+\frac{\Theta^2 T}{\mu L} + \frac{\Theta^2T^2}{L^2}\right)\|\widetilde{B}\|_{L_{[0,T_0]}^{\infty}H_x^k}^4.
$$
Also,
$$|(III)|\leq T\Theta^2\|\widetilde{B}\|_{L_{[0,T_0]}^{\infty}H_x^k}^4.$$
Combining these estimates, we obtain
$$
\|\widetilde{B}(T_0)\|_{H^k}^2\leq \|\widetilde{B}_0\|_{H^k}^2 +
(\|\psi_{+0}\|_{H^k}+\|\psi_{-0}\|_{H^k})\int_0^{T_0}\|\widetilde{B}(t)\|_{H^k}^2 dt
+ \varepsilon\|\widetilde{B}\|_{L_{[0,T_0]}^{\infty}H_x^k}^4,
$$
where $\varepsilon=\tfrac{\Theta^2\mu T^2}{L}+\tfrac{\Theta^2 T}{\mu L}
+ \tfrac{\Theta^2T^2}{L^2}+ \Theta^2 T$. Using our assumptions on $T_0$, $L$
and $\Theta$, we see that $\varepsilon\leq L^{-1/2}\lesssim
e^{-(\|\psi_{+0}\|_{H^k}+\|\psi_{-0}\|_{H^k})T}\|\widetilde{B_0}\|_{H^k}^{-2}$. This fact combined with the
Gronwall inequality and a continuity argument allows us to conclude that
$$\|\widetilde{B}\|_{L_{[0,T_0]}^{\infty}H_x^k}^2\leq 2 e^{(\|\psi_{+0}\|_{H^k}
+\|\psi_{-0}\|_{H^k})T} \|\widetilde{B}_0\|_{H^k}^2.$$
This completes the proof of the proposition.
\end{proof}

Denote by $\widetilde{A}$ the solution of the small-dispersion limit (i.e., $\mu=0$) of~\eqref{e.mzr'}:
\begin{equation}\label{e.step1}
i\partial_x\widetilde{A} = (\psi_{+0}+\psi_{-0})\widetilde{A}
\end{equation}
with $\widetilde{A}(0,x)=\widetilde{B}_0(x)$. In other words,
$$\widetilde{A}(x,t)=e^{-it(\psi_{+0}(x)+\psi_{-0}(x))}\widetilde{B}_0(x).$$

\begin{proposition}\label{p.step3}It holds $\|\widetilde{B}-\widetilde{A}\|_{L_{[0,T_0]}^{\infty}H^k}\lesssim\mu$.
\end{proposition}

\begin{proof}The energy method applied to $\widetilde{B}-\widetilde{A}$ yields
\begin{eqnarray*}
&\|\partial_x^k(\widetilde{B}-\widetilde{A})(T_0)\|_{L^2_x} = \\
&=-2\Re\left(i\mu^2\int_0^{T_0}\int \partial_x^{k+2}\widetilde{B}\overline{\partial_x^k(\widetilde{B}-\widetilde{A})}\right) \\
&-2\Re\left(i\int_0^{T_0}\int \partial_x^{k}(\psi_{\pm0}(x\mp\tfrac{\mu(1\mp c)}{L}t)\widetilde{B}-\psi_{\pm0}(x)\widetilde{A})\overline{\partial_x^k(\widetilde{B}-\widetilde{A})}\right) \\
&-2\Re\left(i\int_0^{T_0}\int \partial_x^{k}(\psi_{\pm}{B})\overline{\partial_x^k(\widetilde{B}-\widetilde{A})}\right) \\
&-2\Re\left(i\Theta^2\int_0^{T_0}\int
\partial_x^k(|\widetilde{B}|^2 \widetilde{B})\overline{\partial_x^k(\widetilde{B}-\widetilde{A})}\right) \\
&= (I) + (II) + (III)+ (IV).
\end{eqnarray*}
Using the Proposition~\ref{p.step2}, we get
$$
|(I)|\leq T^2\mu^4\|\widetilde{B}\|_{L_{[0,T_0]}^{\infty}H_x^{k+2}}^2 + \tfrac{1}{4} \|\partial_x^k(\widetilde{B}-\widetilde{A})\|_{L_{[0,T_0]}^{\infty}L^2_x}^2
\leq cT^2\mu^3 + \tfrac{1}{4} \|\partial_x^k(\widetilde{B}-\widetilde{A})\|_{L_{[0,T_0]}^{\infty}L^2_x}^2,
$$
$$|(III)|\leq T^2\|\widetilde{\psi}_\pm\|_{L_{[0,T_0]}^{\infty}H^k_x}^2
\|\widetilde{B}\|_{L_{[0,T_0]}^{\infty}H^k_x}^2
+ \tfrac{1}{4} \|\partial_x^k(\widetilde{B}-\widetilde{A})\|_{L_{[0,T_0]}^{\infty}L^2_x}^2\leq \tfrac{T^2\Theta^4}{L^2\mu} +
\tfrac{1}{4} \|\partial_x^k(\widetilde{B}-\widetilde{A})\|_{L_{[0,T_0]}^{\infty}L^2_x}^2
$$
and
$$
|(IV)|\leq T^2\Theta^4\|\widetilde{B}\|_{L_{[0,T_0]}^{\infty}H^k_x}^6
+ \tfrac{1}{4} \|\partial_x^k(\widetilde{B}-\widetilde{A})\|_{L_{[0,T_0]}^{\infty}L^2_x}^2
\leq \tfrac{T^2\Theta^4}{\mu^3} +
\tfrac{1}{4} \|\partial_x^k(\widetilde{B}-\widetilde{A})\|_{L_{[0,T_0]}^{\infty}L^2_x}^2.
$$
Furthermore, we can rewrite
\begin{eqnarray*}
&\widetilde{\psi}_{\pm0}(x\mp\tfrac{\mu(1\mp c)}{L}t)\widetilde{B}-\widetilde{\psi}_{\pm0}(x)\widetilde{A} = \\
&= \left(\int_0^{\tfrac{\mp\mu(1\mp c)}{L}t}\partial_x\widetilde{\psi}_{\pm0}(x+t')dt'\right)\widetilde{B} + \widetilde{\psi}_{\pm0}(x)(\widetilde{B}-\widetilde{A}).
\end{eqnarray*}
In particular,
$$|(II)|\leq \tfrac{\mu(1\mp c)T^2}{L}\|\psi_{\pm0}\|_{H^{k+1}}\|\widetilde{B}\|_{L_{[0,T_0]}^{\infty}H^k_x}
\|\widetilde{B}-\widetilde{A}\|_{L_{[0,T_0]}^{\infty}H^k_x} + \|\psi_{\pm0}\|_{H^k}\int_0^{T_0} \|(\widetilde{B}-\widetilde{A})(t)\|_{H^k}^2 dt.$$
Putting these facts together and using Gronwall inequality, we conclude
$$\|\widetilde{B}-\widetilde{A}\|_{L_{[0,T_0]}^{\infty}H^k_x}^2\lesssim e^{cT}\cdot
\left(T^2\mu^3+ \frac{\mu(1\mp c)^2 T^4}{L^2} + \frac{T^2\Theta^4}{L^2\mu} + \frac{T^2\Theta^4}{\mu^3}\right).$$
This ends the proof of the proposition.
\end{proof}

Define
$$
B(x,t):=L\Theta e^{-ic^2t}e^{icx}\widetilde{B}(L\mu(x-ct),L^2t),
$$
$$
\psi_{\pm}(x,t) = L^2\widetilde{\psi}_{\pm0}(L\mu(x\mp t)) + L^2\widetilde{\psi}_{\pm}(L\mu(x-ct),L^2t),
$$
where $(\widetilde{B},\widetilde{\psi}_\pm)$ is a solution of~\eqref{e.mzr'}. It is quite straightforward
to check that $(B,\psi_{\pm})$
solves the Zakharov-Rubenchik system~\eqref{e.zr'} with initial data
$B(x,0)=L\Theta e^{icx} \widetilde{B}_0(L\mu x)$, $\psi_{\pm}(x,0)=L^2\widetilde{\psi}_{\pm0}(L\mu x)$.

Since $\|\widetilde{B}(x,t)\|_{L^2_x} = \|\widetilde{B_0}\|_{L^2_x}$ for all $t$, we get
$$\|B(x,t)\|_{L^2_x} = \frac{L^{1/2}\Theta}{\mu^{1/2}}\|\widetilde{B_0}\|_{L^2_x}.$$
Also, if $\widehat{\psi_{+0}}(\xi)=0$ for $|\xi|\leq 1$ and $\mu L\geq 1$, we see
$$\|\psi_{+}(x,0)\|_{H^l_x}\leq L^{\tfrac{3}{2}+l} \mu^{l-\tfrac{1}{2}}\|\psi_{+0}\|_{H^l_x}.$$
Therefore, if $l<-3/2$ and $L\geq\mu^{-\alpha}$, where $\alpha=\max\{5,\tfrac{1/2-l}{-l-3/2}\}$, we obtain
$$\|\psi_{+}(x,0)\|_{H^l_x}\leq \|\widetilde{\psi}_{+0}\|_{H^l_x}.$$

At this point, we are ready to complete the proof of Theorem~\ref{t.ill-decoherence}. We fix $M\gg 1$ and
$0<\mu\ll 1$ to be chosen later and we put $T=|\log\mu|\cdot M^{-2}$,
$$L_1 = M, \quad L_2 = \sqrt{\frac{\pi}{2T}+ M^2}.$$
Observe that $e^{iT(L_2^2-L_1^2)}=i$ and $L_2/L_1 = \sqrt{\tfrac{\pi}{2|\log\mu|}+1}\to 1$ (uniformly on $M$).
Take $\Theta^2 = \mu \cdot M^{-1}$. Fix $\widetilde{B}_0$ a Schwartz function such that $\widetilde{B}_0(x)=1$
for $|x|\leq 1$, and $\widetilde{\psi}_{+0}(x) = \cos (3x)\sin(x)/x$ and $\widetilde{\psi}_{-0}(x)=0$.

Next we consider $\widetilde{B}_1$ and $\widetilde{B}_2$ solutions of~\eqref{e.mzr'} with parameters $L_1$
and $L_2$ (resp.) and same initial data $\widetilde{B}_0$, $\widetilde{\psi}_{\pm0}$. From the previous
discussion, for $j=1,2$,
$$
\|B_j(x,t)\|_{L^2_x} \sim 1, \quad \|\psi_{+,j}(x,t)\|_{H^l_x}\leq 1,
\quad \|\psi_{-,j}(x,t)\|_{H^s_x}\leq 1,
$$
whenever $M\geq \mu^{-\alpha}$.

On the other hand, using the definition of $\Theta$,
$$\|B_2(x,t)-B_1(x,t)\|_{L^2_x} = \|\tfrac{L_2}{L_1}\widetilde{B}_2(\tfrac{L_2}{L_1}x,L_2^2t) - \widetilde{B}_1(x,L_1^2t)\|_{L^2_x}.$$
Since $L_2/L_1\to 1$ uniformly on $M$ and $\|\widetilde{B}_j(x,t)\|_{L^2_x} = \|\widetilde{B}_0\|_{L^2_x}$ for all $t$, we can select $\mu=\mu(\delta)>0$ sufficiently small so that
$$\|B_2(x,t)-B_1(x,t)\|_{L^2_x} = \|\widetilde{B}_2(x,L_2^2t) - \widetilde{B}_1(x,L_1^2t)\|_{L^2_x} + O(\delta).$$
Using the proposition~\ref{p.step3}, it follows that
$$\|B_2(x,t)-B_1(x,t)\|_{L^2_x} = \|(e^{i(L_2^2-L_1^2)t\widetilde{\psi}_{+0}(x)}-1)\widetilde{B}_0\|_{L^2_x} + O(\delta).$$
Hence, $\|B_2(x,0)-B_1(x,0)\|_{L^2_x} \lesssim \delta$ but $\|B_2(x,T)-B_1(x,T)\|_{L^2_x} \sim 1$ (because $e^{iT(L_2^2-L_1^2)}=i$).
Since $T=|\log\mu|\cdot M^{-2}\leq |\log\mu|\cdot\mu^{10}\to 0$ as $\mu\to 0$, the proof of the
Theorem~\ref{t.ill-decoherence} is complete.

\end{document}